\documentclass[11pt,a4paper,leqno,twoside, headinclude]{amsart}

\usepackage[tracking=true]{microtype}
\DeclareMicrotypeSet*[tracking]{my}%
  { font = */*/*/sc/* }%
\SetTracking{ encoding = *, shape = sc }{ 45 }

\title{Special Holonomy on Special Spaces}
\def\titl{Special Holonomy on Special Spaces}
\def\auth{Manuel Amann}
\date{March 4th, 2013}

\thanks{}

\author{\auth}

\usepackage[english]{babel}

\usepackage[colorlinks,pdftex, plainpages=false]{hyperref}
\hypersetup{pdftitle=\titl, pdfauthor=\auth, pdftoolbar=false,
plainpages=false, hyperindex=true, pdfdisplaydoctitle=true}

\hypersetup{colorlinks,%
citecolor=black,%
filecolor=black,%
linkcolor=black,%
urlcolor=black,%
pdftex}

\usepackage{color}

\usepackage{amsmath, amssymb, amscd, amsthm}
\usepackage{stmaryrd}
\usepackage{latexsym}
\usepackage[all]{xy}
\usepackage{pb-diagram}
\usepackage{rotating}
\usepackage{multicol}
\usepackage{lscape}
\usepackage{wasysym}
\usepackage{longtable}
\usepackage{enumerate}

\xyoption{all}

\newtheorem{theo}{Theorem}[section]
\newtheorem{main}{Theorem}

\newtheorem*{main*}{Theorem}

\newtheorem*{mainprop*}{Proposition}

\newtheorem{mainconj}{Conjecture}

\newtheorem{prop}[theo]{Proposition}
\newtheorem{defi2}[theo]{Definition}
\newtheorem*{defi2*}{Definition}

\newenvironment{defi*}{\begin{defi2*}\normalfont}{\end{defi2*}}

\newenvironment{defin*}[1]{\begin{defi2*}[#1]\normalfont}{\end{defi2*}}
\newtheorem*{rem2*}{Remark}
\newenvironment{rem*}{\begin{rem2*}\normalfont}{\hfill$\boxbox$\end{rem2*}}
\newtheorem{rem2}[theo]{Remark}
\newenvironment{rem}{\begin{rem2}\normalfont}{\hfill$\boxbox$\end{rem2}}

\newtheorem{lemma}[theo]{Lemma}
\newtheorem{cor}[theo]{Corollary}
\newtheorem*{cor*}{Corollary}

\newtheorem*{conj*}{Conjecture}
\newtheorem*{theo*}{Theorem}
\newtheorem*{ques*}{Question}
\newtheorem*{mi2}{Main Idea}

\newtheorem{ex2}[theo]{Example}

\newtheorem{exer2}[theo]{Exercise}

\newtheorem{alg2}[theo]{Algorithm}

\newcommand{\cc}{{\mathbb{C}}}                                     
\newcommand{\kk}{{\mathbb{K}}}                                     
\newcommand{\hh}{{\mathbb{H}}}                                     
\newcommand{\qq}{{\mathbb{Q}}}                                     
\newcommand{\rr}{{\mathbb{R}}}                                     
\newcommand{\pp}{{\mathbf{P}}}                                     
\newcommand{\Gr}{{\mathbf{Gr}}}                                    
\newcommand{\s}{{\mathbb{S}}}                                      
\newcommand{\zz}{{\mathbb{Z}}}                                     
\newcommand{\GL}{{\mathbf{GL}}}                                    
\newcommand{\SO}{{\mathbf{SO}}}                                    
\newcommand{\U}{{\mathbf{U}}}                                      
\newcommand{\SU}{{\mathbf{SU}}}                                    
\newcommand{\Sp}{{\mathbf{Sp}}}                                    
\newcommand{\E}{{\mathbf{E}}}                                      
\newcommand{\F}{{\mathbf{F}}}                                      
\newcommand{\G}{{\mathbf{G}}}                                      
\newcommand{\Spin}{{\mathbf{Spin}}}                                
\newcommand{\dif} {{\operatorname{d}}}                             
\newcommand{\In} {{\,\subseteq\,}}                                 
\newcommand{\Ni} {\,{\supseteq}\,}                                 
\newcommand{\Aut}{{\operatorname{Aut}}}                            
\newcommand{\Hol}{{\operatorname{Hol}}}                            
\newcommand{\id}{{\operatorname{id}}}                              
\newcommand{\ADR}{{\operatorname{A_{DR}}}}                         
\newcommand{\Sym}{{\operatorname{Sym}}}                            
\newcommand{\cat} {{\operatorname{cat}}}                           
\newcommand{\rk}{{\operatorname{rk\,}}}                            
\newcommand{\diag}{{\operatorname{diag}}}                          
\newcommand{\co}{\colon\thinspace}                                 

\newcommand{\comment}[1]{}                                         
\newcommand{\xto}[1]{\xrightarrow{#1}}                             
\newcommand{\hto}[1]{\overset{#1}{\hookrightarrow}}                

\newcommand{\biq}[2]{#1\;\!\!\!\sslash \;\!\!\!#2}                 
\newcommand{\str}{\noindent\textbf{Structure of the article. }}    
\newcommand{\odd}{\textrm{odd}}                                    
\newcommand{\even}{\textrm{even}}                                  

\newenvironment{prf}{\begin{proof}[\textsc{Proof}]} {\end{proof}}     

\begin{document}

\maketitle \thispagestyle{empty}


\begin{abstract}
We characterise simply-connected biquotients which potentially admit metrics of holonomy $\G_2$.
We prove that there are at most three real homotopy types of rationally elliptic such manifolds---all of them being formal. In the course of this examination we classify rationally elliptic homotopy types and characterise $7$-dimensional simply-connected biquotients from a rational point of view.
Moreover, we also investigate further manifolds of special holonomy, like manifolds of holonomy $\Spin(7)$ or $\Sp(n)\Sp(1)$ in special situations provided by rational ellipticity or geometric formality.
\end{abstract}


\section*{Introduction}

Due to Berger et al.~simply-connected non-symmetric irreducible Riemannian manifolds $(M,g)$ fall into a very short list of possible holonomy types. That is
\begin{align*}
\Hol(M,g)\in \{\SO(n), \U(n), \SU(n), \Sp(n), \Sp(n)\Sp(1), \G_2, \Spin(7)\}
\end{align*}
where $n=\dim M$, $n=\dim M/2$ and $n=\dim M/4$ respectively. (For the definition of the holonomy group see Section \ref{sec02}.) Ignoring the ``generic case'' of $\SO(n)$-holonomy, one obtains manifolds of ``special holonomy'' which, interestingly, reveal several pretty special topological features. One characteristic all these manifolds share is a Lefschetz-like property which (formulated differently and more or less strong in the respective cases) underlies the structure of Poincar\'e duality of their real cohomology algebras.

This feature has a ``formalising tendency''. Indeed, another interesting conjectured property (cf.~for example \cite{AK12}) is the formality of simply-connected manifolds of special holonomy. In simplified terms, this means that their rational homotopy type is (up to the application of an algorithm) ``the same'' as their rational cohomology algebra.

This was verified for K\"ahler manifolds in \cite{DGMS75} and for positive quaternion K\"ahler manifolds in \cite{AK12}. Since $\SU(n)\In \U(n)$ and $\Sp(n)\In \U(2n)$ this basically leaves open the cases of $\G_2$-manifolds and $\Spin(7)$-manifolds. In \cite{Cav05} it is shown that the known topological properties of simply-connected $\G_2$-manifolds are not enough to prove their formality. In this article we shall prove this property under the additional assumption of rational ellipticity, i.e.~assuming that only finitely many homotopy groups are not entirely torsion.

This property, however, just serves as a first motivation for our investigation, since we obtain a much more general result: We classify the real homotopy type of elliptic $\G_2$-manifolds in the rationally elliptic case and we show that there are no such manifolds of holonomy $\Spin(7)$. Actually, for this we provide a real classification of \emph{all} $7$-dimensional spaces. Manifolds of holonomy $\G_2$ and $\Spin(7)$ were investigated in \cite{Joy00} and both the methods of investigation and the results reveal some similarities between them. So we try to deal with them at the same time.

\begin{main}\label{theo03}
\begin{itemize}
\item
A rationally elliptic
compact manifold of dimension $7$ does not admit a metric $g$ with
$\Hol_g(M)=\Spin(7)$.
\item
Neither does any simply-connected eight-dimensional rationally elliptic manifold admit a metric with $\Hol_g(M)= \SU(4)$. No $8$-dimensional compact homogeneous space admits a metric of holonomy $\Sp(2)$.
\item There are at most three different real types (the formal types $\s^4\times \s^3$, $\cc\pp^2\times \s^3$ and $\cc\pp^2\#\cc\pp^2\times \s^3$) of simply-connected rationally elliptic manifolds admitting a metric of holonomy $\G_2$.
\item Every simply-connected irreducible rationally elliptic $8$-manifold admitting a metric of holonomy contained in $\Spin(7)$ is formal.
\end{itemize}
\end{main}
As mentioned above, part of this theorem is proven via a stronger result, namely the \emph{classification of real homotopy types of simply-connected rationally elliptic spaces}. We show that there are only the types
\begin{align*}
&\s^7, \s^4\times \s^3, \s^2\times \s^5, \cc\pp^2\times \s^3, \s^2\times\s^2\times\s^3, \\&\s^3\times \cc\pp^2\#\cc\pp^2,
\s^3\tilde\times (\s^2\times\s^2)
\end{align*}
where $\s^3\tilde\times (\s^2\times\s^2)$ is the only non-trivial real/rational $\s^3$-bundle over $\s^2\times \s^2$---see Lemma \ref{lemma02}.

A large and prominent class of rationally elliptic spaces are biquotients---we recall the definition in Section \ref{sec01}. Thus it seems natural to sort out those biquotients which might admit a metric of $\G_2$-holonomy.
\begin{main}\label{theo01}
Let $M$ be a simply-connected biquotient admitting a metric of holonomy $\G_2$.
Then $M$ is $\biq{\Sp(1)\times \Sp(1)\times \Sp(1)}{\s^1\times \s^1}$. In particular, there is no such biquotient of a simple Lie group.

There is no homogeneous space admitting a metric of holonomy $\G_2$.
\end{main}

This result is a consequence of the subsequent classification result depicted in Theorem \ref{theo02}. We point the reader to the Ph.D.~Thesis \cite{Dev11} by Jason DeVito, where biquotients up to dimension $7$ were classified and where further properties of these spaces were investigated. Several examinations there go deeper than what we need in this article; in particular, concrete actions of the denominator group are studied. Nonetheless, for the convenience of the reader and since the classification focus in \cite{Dev11} is clearly not a real one, we produce a purely real classification result from scratch following an approach by Totaro (cf.~\cite{Tot02}) and are convinced that this approach will be easier to follow than an adaptation of the results in \cite{Dev11}.
\begin{main}\label{theo02}
A simply-connected seven-dimensional biquotient is as in Table \ref{table07}.
\end{main}
\begin{table}[h]
\centering \caption{biquotients $\biq{G}{H}$ in dimension $7$}
\label{table07}
\begin{tabular}{c@{\hspace{4mm}}| @{\hspace{4mm}}c@{\hspace{4mm}}}
real type & $\biq{G}{H}$ \\
\hline
$\s^7$ & $\biq{\SU(4)}{\SU(3)}$, $\biq{\Sp(2)}{\Sp(1)}$,\\& $\Spin(7)/\G_2, \SO(8)/\SO(7)$ \\
$\s^4\times \s^3$ & $(\Sp(2)/\Sp(1)) \times \Sp(1)$\\
$\s^2 \times \s^5$ &  $\biq{\SU(3)}{\s^1}$, $\biq{\SU(3)\times \Sp(1)}{\s^1\times \Sp(1)}$,\\&  $\biq{\SU(4)\times \Sp(1)}{\Sp(2)\times \s^1}$\\
$\cc\pp^2 \times \s^3$ &$\biq{\SU(3)\times \Sp(1)}{\s^1\times \Sp(1)}$\\
$\s^2\times\s^2\times\s^3$ &  $\biq{\Sp(1)\times \Sp(1)\times \Sp(1)}{\s^1\times \s^1}$\\
$\s^3\times \cc\pp^2\#\cc\pp^2$ &  $\biq{\Sp(1)\times \Sp(1)\times \Sp(1)}{\s^1\times \s^1}$\\
$\s^3\tilde\times (\s^2\times\s^2)$ &  $\biq{\Sp(1)\times \Sp(1)\times \Sp(1)}{\s^1\times \s^1}$
\end{tabular}
\end{table}

We point the reader to \cite{Tot03} where infinitely many rationally distinct biquotients in dimension six, in particular, are constructed.

\vspace{3mm}

Moreover, we provide several classification results for positive quaternion K\"ahler manifolds and we also investigate special holonomy in combination with geometric formality---see Sections \ref{sec04} and  \ref{sec05}.

For example (see Theorem \ref{theo04}) we show that every positive quaternion K\"ahler manifold diffeomorphic to a biquotient of the rational homotopy type of a compact rank one symmetric space is homothetic to
\begin{align*}
\hh\pp^n, \qquad \G_2/\SO(4) \quad \textrm{or} \quad \Gr_2(\cc^4)
\end{align*}

Let us state the following conjecture (which follows easily from the classical LeBrun--Salamon conjecture---cf.~Section \ref{sec04}---in this context) and which spurs this investigation.
\begin{conj*}
A positive quaternion K\"ahler manifold is (rationally) elliptic and geometrically formal.
\end{conj*}
We shall also motivate this conjecture at the end of Section \ref{sec08}. For a definition of \emph{ellipticity} (over $\zz$ instead of $\qq$) we point the reader to \cite{FHT91}.

\vspace{3mm}

We end this article by a comment on a change of the coefficient field---see Theorem \ref{theo05}. It implies, in particular, that if in a certain dimension there are infinitely many complex homotopy types of compact manifolds, then there are also already infinitely many real homotopy types.

\vspace{3mm}

\str In Section \ref{sec01} we recall some definitions and some well-known results from special holonomy. In Section \ref{sec02} we provide a classification of simply-connected $7$-dimensional real homotopy types, which finally yields a proof of Theorem \ref{theo03}. In Section \ref{sec03} we prove Theorem \ref{theo01}, the real classification of $7$-dimensional simply-connected biquotients. From this result Theorem \ref{theo02} then follows via some further arguments. In Section \ref{sec04} we deal with the properties of rationally elliptic positive quaternion K\"ahler manifolds before we investigate geometric formality in the context of special holonomy in Section \ref{sec05}. In Section \ref{sec09} we provide the result on the change of the coefficient field.

\vspace{3mm}


\section{Preliminaries} \label{sec01}

Recall the definition of the \emph{holonomy group} of a Riemannian manifolds $(M,g)$ as the group
\begin{align*}
\Hol_x(M,g)=\{P_\gamma\mid \gamma \textrm{ is a closed loop based at
}x\} \In \GL(T_xM)
\end{align*}
Here $P_\gamma$ denotes parallel transport and the holonomy group $\Hol_g(M):=\Hol_x(M,g)$ is independent of the chosen base point $x$ up to inner automorphism.

As indicated in the introduction Berger's theorem allows to speak of manifolds of ``special holonomy'' $\U(n),\SU(n),\Sp(n),\Sp(n)\Sp(1),\G_2,\Spin(7)$, and these manifolds bear remarkable topological features. Most prominently, they share Lefschetz-like properties. We shall illustrate this in the case of manifolds of holonomy $\G_2$ or $\Spin(7)$.

A compact manifold $M$ of holonomy $\G_2$ is orientable, spin, $\pi_1(M)$ is finite and the first Pontryagin class does not vanish, i.e.~$p_1(M)\neq 0$---see \cite[Theorem 10.2.8, p.~247]{Joy00}.
In \cite[Proposition 10.2.6, p.~246]{Joy00} it is proven that a
compact manifold $M$ with
$\Hol(M)=\G_2$, satisfies
\begin{align}\label{PUReqn09}
\langle a\cup a\cup[\omega],[M]\rangle <0
\end{align}
for every non-zero $a\in H^2(M;\rr)$ and with respect to the $3$-form
$\omega$ defining the $\G_2$-structure (cf.~\cite[Definition 10.1.1, p.~242]{Joy00}).

The analogue holds for compact $\Spin(7)$-manifolds
(cf.~\cite[Proposition 10.6.6, p.~261]{Joy00}) with the respective $4$-form $\omega$
(cf.~\cite[Definition 10.5.1, p.~255]{Joy00}). So, in particular, combining this
with Poincar\'e duality one obtains Lefschetz-like properties for
manifolds with holonomy $\G_2$ or $\Spin(7)$. It holds
\begin{align*}
\cup \omega: H^2(M;\rr)\xto{\cong} H^{n-2}(M;\rr)
\end{align*}
for the closed $3$-form $\omega$ defining the $\G_2$-structure and
$n=7$; respectively the $4$-form $\omega$ defining the
$\Spin(7)$-structure and $n=8$.

If $\Hol(M)=\Spin(7)$, \cite[Theorem 10.6.8, p.~261]{Joy00} gives
further strong topological restrictions; there is the following
relation on Betti numbers:
\begin{align}\label{PUReqn08}
 b_3+b_4^+=b_2+2 b_4^-+25
\end{align}
where $b_4^+ +b_4^-=b_4$ and $b_4^+\geq 1$. The manifold is
simply-connected and spin with $\hat
A(M)[M]=1$---cf.~\cite[Equation 10.25, p.~260]{Joy00} together
with \cite[Proposition 10.6.5, p.~260]{Joy00}. A direct consequence of the latter
result is that there are no smooth effective $\s^1$-actions upon $M$
due to \cite{AH70}. In particular, we shall not find any homogeneous
space that admits a metric with such
holonomy.

In the case of holonomy equal to $\SU(4)$ we cite from \cite[Theorem p.~113]{Sal89} that
\begin{align}\label{PUReqn10}
b_3+b_4^+\geq 50
\end{align}

\vspace{3mm}

Finally, we recall the definition of a biquotient:
Let $G$ be a compact connected Lie group and let
$H\In G\times G$ be a closed (Lie) subgroup.

Then $H$ acts on $G$ on
the left by $(h_1,h_2)\cdot g=h_1gh_2^{-1}$. The orbit space $G/H$ of this
action is called the \emph{biquotient} $\biq{G}{H}$ of $G$ by $H$. If the action of $H$ on $G$ is \emph{free}, then
$\biq{G}{H}$ possesses a manifold structure. This is the only case
we shall consider.

\section{Proof of Theorem \ref{theo03}}\label{sec02}

In this section we provide a proof of Theorem \ref{theo03}. An essential ingredient in this proof will be a more general result, namely the real classification of $7$-dimensional rationally elliptic spaces.

Since the groups $\pi_*(M)\otimes \rr$ vanish if and only if so do the $\pi_*(M)\otimes \qq$ ``real ellipticity'' is equivalent to rational ellipticity. However, let us use the latter, more common terminology.

A rationally elliptic space satisfies strong restrictions on the configuration of both its rational homotopy groups and Betti numbers. We recall the homotopy Euler characteristic
\begin{align*}
\chi_\pi(M)=\dim \pi_\odd(M)\otimes \qq-\dim \pi_\even(M)\otimes \qq
\end{align*}

We cite the following relations from \cite[p.~434]{FHT01}.
\begin{align}
\label{eqn01}
\sum_{i} \deg x_i &\leq 2\cdot \dim M -1
\\ \label{eqn02}
\sum_{i} \deg y_i &\leq \dim M \\
\label{eqn03} \dim \pi_*(M)\otimes \qq &\leq \dim M\\
\label{eqn04} \sum_i \deg x_i-\sum_j (\deg y_j -1)&=\dim M\\
\label{eqn05} \chi_\pi M\geq 0 \textrm{ and } \chi M\geq 0 \textrm{ and } \chi_\pi M>0 &\Leftrightarrow \chi M=0
\end{align}
Here the $x_i$ form a homogeneous basis of
$\pi_\textrm{odd}\otimes \qq$ and the $y_i$ form a homogeneous basis
of $\pi_\textrm{even}\otimes \qq$.

These relations permit to prove the following classification result. Compare this to the classification (see \cite[Example 3.8, p.~108]{FOT08}, \cite[Lemma 3.2, p.~426]{PP03}) of $4$-manifolds, which fall into the rational homotopy types
\begin{align*}
\s^4, \cc\pp^2, \s^2\times \s^2,  \cc\pp^2\#\cc\pp^2
\end{align*}
and the homeomorphism types
\begin{align*}
\s^4, \cc\pp^2, \s^2\times \s^2, \cc\pp^2\#\cc\pp^2, \cc\pp^2\#\overline{\cc\pp^2}
\end{align*}
Note that $\cc\pp^2\#\overline{\cc\pp^2}\simeq_\qq \s^2\times \s^2$.
\begin{prop}\label{lemma02}
\begin{enumerate}
\item
A simply-connected rationally elliptic $4$-dim\-en\-sional space $M$ has the real homotopy type of $\s^4$, $\cc\pp^2$, $\s^2\times \s^2$, $\cc\pp^2\#\cc\pp^2$.
\item
The real homotopy type determines the rational homotopy type in the case of $\s^4$, $\cc\pp^2$.
\item
There are infinitely many rational homotopy types realising the real homotopy types $\cc\pp^2\#\cc\pp^2$ and $\cc\pp^2\#\overline{\cc\pp^2}\simeq_\qq \s^2\times \s^2$ respectively. So there are infinitely many different rational homotopy types of orbifolds in each of these cases.
\item There are infinitely many rational homotopy types of simply-connected smooth compact simply-connected seven-manifolds.
\item A simply-connected $7$-dimensional rationally elliptic space $M$ has one of the following real homotopy types
\begin{align*}
&\s^7, \s^4\times \s^3, \s^2\times \s^5, \cc\pp^2\times \s^3, \s^2\times\s^2\times\s^3, \\&\s^3\times \cc\pp^2\#\cc\pp^2,
\s^3\tilde\times (\s^2\times\s^2)
\end{align*}
where $\s^3\tilde\times (\s^2\times\s^2)$ is the only non-trivial rational/real $\s^3$-bundle over $\s^2\times \s^2$.
\end{enumerate}
\end{prop}

\begin{prf}

\textbf{ad (1),(2)}
We can easily compute that the rational homotopy groups of $M$ fall into one of the following categories
\begin{enumerate}
\item $\pi_*(M)\otimes \qq=\pi_4(M)\otimes \qq \oplus \pi_7(M)\otimes \qq$ with $\dim \pi_4(M)\otimes \qq=\dim\pi_7(M)\otimes \qq=1$.
\item $\pi_*(M)\otimes \qq=\pi_2(M)\otimes \qq \oplus \pi_5(M)\otimes \qq$ with $\dim \pi_2(M)\otimes \qq=\dim \pi_5(M)\otimes \qq=1$.
\item $\pi_*(M)\otimes \qq=\pi_2(M)\otimes \qq \oplus \pi_3(M)\otimes \qq$ with $\dim \pi_2(M)\otimes \qq=\dim\pi_3(M)\otimes \qq=2$.
\end{enumerate}
Indeed, since a rationally elliptic space satisfies Poincar\'e duality, we obtain $H^1(M)=H^3(M)=0$ and $\chi(M)>0$. This implies that $\chi_\pi M=0$ due to \eqref{eqn05}. According to \eqref{eqn01} we have that $\dim \pi_\odd(M)\leq 2$ and \eqref{eqn01} and \eqref{eqn04} then yield the potential configurations: That is, if $\pi_7(M)\otimes \qq\neq 0$, then is $\qq$ and $\pi_\odd(M)\otimes \qq=\pi_7(M)\otimes \qq$. The same holds for $\pi_5(M)\otimes \qq$. The fact that $H^i(M)=0$ for $i>4$ then implies that $H^*(M)=H^4(M)=\qq$ in the first case and $H^*(M)=\qq[x]/x^3$ for $\deg x=2$ in the second one. In the remaining case, Case (3), we necessarily have $\dim\pi_\odd(M)\otimes \qq=\dim \pi_3(M)\otimes \qq=2$ and $\dim \pi_\even(M)\otimes \qq=\dim \pi_2(M)\otimes \qq=2$ by \eqref{eqn04}. Positive Euler characteristic then implies that $H^*(M)$ is a truncated polynomial ring---by a regular sequence---generated by two elements in degree $2$.

It is now easy to see that Case (1) implies that $M$ has the rational type of $\s^4$. Case (2) can only be realised by $\cc\pp^2$.

For Case $(3)$ we form the minimal model of such a space as $(\Lambda V,\dif)$ with $V=\langle a,b,x,y\rangle$, $\deg a=\deg b=2$, $\deg x=\deg y=3$ and $\dif$ vanishes in degree $2$ and is injective on degree $3$.

Via a change of basis and the fact that the differentials of $x$ and $y$ must form a regular sequence, Case (3) now falls into several subcases (up to isomorphism).
\begin{enumerate}
\item[(3.1)] $\dif(x)=a^2+sb^2$, $\dif(y)=ab$
\item[(3.2)] $\dif(x)=a^2+sab$, $\dif(y)=b^2$
\item[(3.3)] $\dif(x)=a^2+sab$, $\dif(y)=b^2+tab$
\end{enumerate}
with $s\neq 0$ in Case (3.1), $s,t\neq 0$ in Case (3.3).

Case (3.1) has the real homotopy type of $\cc\pp^2\#\cc\pp^2$ via the isomorphism $a\mapsto a$, $b \mapsto (1/\sqrt{s}) b$, $x\mapsto x$, $y\mapsto (1/\sqrt{s}) y$ if $s>0$, and the real homotopy type of $\cc\pp^2\#\overline{\cc\pp^2}$ via the isomorphism $a\mapsto a$, $b \mapsto (1/\sqrt{-s}) b$, $x\mapsto x$, $y\mapsto (1/\sqrt{-s}) y$ if $s<0$.

Case (3.2) corresponds to $\s^2_{a+(s/2) b}\times \s^2_b$ up to the isomorphism on cohomology induced by $a\mapsto a+\frac{s}{2}b$, $b\mapsto b$. This is an isomorphism of \emph{rational} homotopy types. The space $\s^2_a\times \s^2_b$ is rationally equivalent to $\cc\pp^2\#\overline{\cc\pp^2}$ as the isomorphism induced on minimal models by $a\mapsto a+b, b\mapsto a-b$ shows.

As for Case (3.3) we differentiate between the following subcases
\begin{itemize}
\item[(3.3.1)] $1-st\geq 0$
\item[(3.3.2)] $1-st<0$
\end{itemize}
In Case (3.3.1) we compute that we have the equivalent relations $a^2+(s+kt)ab+kb^2$ for $k\in \qq$ and $b^2+tab=0$.
Obviously, $t=0$ yields the previous case. Since $1-st\geq 0$, the equations $v^2=k$ and $2v=s+kt$ yield the quadratic equation $tv^2-2v+s=0$, which is solvable in $v\in \rr$ via $v=\frac{1\pm \sqrt{1-st}}{t}$. In other words, we may replace the relations by yet further equivalent ones: $(a+vb)^2=0$, $b^2+tab=0$. Setting $c:=a+vb$ we recognise Case (3.2) in the relations
\begin{align*}
c^2=0, \qquad (1-v)b^2+tbc=0.
\end{align*}
unless $1-v= 0 \Leftrightarrow s=1/t$. If $st=1$, then the relations yield $ac=0$ and $c^2=0$, which is a contradiction to a regular sequence.
The setting then is equivalent to $\s^2\times\s^2$.

For Case (3.3.2) we specify the isomorphism
\begin{align*}
\varphi\co H(\Lambda V,\dif)\to H^*(\cc\pp^2\#\cc\pp^2)
\end{align*}
defined by $a\mapsto a+b$, $b\mapsto -\frac{1+\sqrt{st-1}}{s}a  + \frac{-1+\sqrt{st-1}}{s}b$. That is, it is represented by the matrix
\begin{align*}
\begin{pmatrix}
1 & 1 \\
-\frac{1+\sqrt{st-1}}{s} & \frac{-1+\sqrt{st-1}}{s}
\end{pmatrix}
\end{align*}
with determinant $\frac{2\sqrt{st-1}}{s}\neq 0$, since $1-st<0$.
Note that the existence of the isomorphism $\varphi$ is equivalent to specifying an isomorphism of rational types. For this we note that both spaces are hyperformal, i.e.~in particular, intrinsically formal and there is exactly one rational homotopy type realising the respective rational cohomology algebras.

The morphism $\varphi$ is a well-defined morphism of algebras: We compute that
\begin{align*}
0&=\varphi(a^2+sab)
\\&=\varphi(a)^2+s\varphi(a)\varphi(b)
\\&=(a+b)^2+s(a+b)\bigg(-\frac{1+\sqrt{st-1}}{s}a  + \frac{-1+\sqrt{st-1}}{s}b\bigg)
\\&=2a^2-s\cdot (2/s)a^2
\\&=0
\end{align*}
and analogously for the second relation $0=\varphi(b^2+tab)=\varphi(b)^2+t\varphi(a)\varphi(b)=(2t/s)a^2-(2t/s) a^2=0$. This proves that Case (3.3.2) is equivalent to $\cc\pp^2\#\cc\pp^2$.

\vspace{5mm}

\textbf{ad (3)}
Let us now prove that there are infinitely many rational types realising the real type of $\cc\pp^2\#\cc\pp^2$ respectively of $\cc\pp^2\#\overline{\cc\pp^2}\simeq_\qq \s^2\times \s^2$.
We make the following observation: We can define an isomorphism $\psi_m=(m\cdot)$ of cochain algebras on $(V,\dif)$ via $a\mapsto ma, b\mapsto mb, x\mapsto m^2x, y\mapsto m^2y$ for $m\in \qq\setminus \{0\}$. Indeed, we obtain
\begin{align*}
\psi_m(\dif x)&=\psi(k_1 a^2+k_2 ab+k_3b^2)
\\&=m^2(k_1a^2+k_2ab+k_3b^2)
\\&=\dif (m^2x)
\\&=\dif(\psi_m(x))
\end{align*}
and analogously for $y$.

Let $(\Lambda V,\dif_1)$ and $(\Lambda V,\dif_2)$ be two minimal Sullivan algebras (over $\qq)$ of spaces $M_1$ and $M_2$ within Case (3.1). The first relation for $(\Lambda V,\dif_1)$ is given by $a^2+sb^2$, the one for $(\Lambda V,\dif_2)$ by $a^2+tb^2$.
Now suppose that $\varphi\co H^*(M_1)\to H^*(M_2)$ is an isomorphism between them defined by $a\mapsto k_1a+k_2b$, $b\mapsto k_3a+k_4b$ with $k_1,k_2,k_3,k_4\in \qq$.  (By abuse of notation we do not differentiate between the morphism in cohomology and the one on Sullivan algebras.)
Using the composition $\psi_{1/k_1}\circ\varphi$ with the automorphism $\psi_{1/k_1}$ we may assume that $k_1=1$ unless $k_1=0$.

If $k_1=0$, we apply the same argument to restrict to the case when $k_2=1$---$\varphi(a)$ may not vanish entirely then. Now we have $\varphi(a)=b$, $\varphi(b)=k_3a+k_4 b$. Since $\varphi$ is multiplicative, we compute $0=\varphi(ab)=b(k_3a+k_4b)=k_4b^2$ and $0=\varphi(a^2+sb^2)=b^2+sk_3^2a^2=(1-tsk_3^2)b^2$. The relation $1-tsk_3^2=0$ is equivalent to
\begin{align}\label{eqn08}
k_3=\pm \sqrt{1/(st)}
\end{align}

Let us now deal with the case when $k_1=1$.
Let $(\Lambda V,\dif_1)$ be given by $\dif_1 x=a^2+sb^2$, $\dif_1 y=ab$, and $(\Lambda V,\dif_2)$ by $\dif_2 x=a^2+tb^2$, $\dif_2 y=ab$. The fact that $\varphi$ is multiplicative yields that $0=\varphi(ab)=(a+k_2b)(k_3a+k_4b)$ which is equivalent to
\begin{align}\label{eqn06}
-tk_3+k_2k_4=0
\end{align}
Analogously, we derive that $(a+k_2b)^2+s(k_3+k_4b)^2=0$, which is equivalent to
\begin{align}\label{eqn07}
-t+k_2^2+s(-tk_3^2+k_4^2)=0
\end{align}
Equation \eqref{eqn06} yields $k_3=k_2k_4/t$. (Obviously, $s,t\neq 0$ for the relations to form a regular sequence.) Using equation \eqref{eqn07} we derive $k_2=\pm\sqrt{t}$ and $k_3=\pm \frac{k_4}{\sqrt{t}}$. The determinant $\det \varphi$ then computes as $1\cdot k_4-\frac{\pm k_4}{\sqrt{t}}\cdot (\pm\sqrt{t})=0$ and $\varphi$ is not invertible; a contradiction.

This means that an isomorphism is subjected to the relations leading to Equation \eqref{eqn08}. Now let the $s\neq t$ run over the prime numbers. Then $k_3$ is never rational and $M_1$ and $M_2$ cannot have the same rational homotopy type. More precisely, we can find infinitely many parameters $s,t$ belonging to $M_1, M_2$ such that these spaces have the same real homotopy type, but are not equivalent over the rationals.

From this equation we can also easily see that $\cc\pp^2\#\cc\pp^2$ and $\cc\pp^2\#\overline{\cc\pp^2}$ are not even equivalent over the reals.
 
Every rational type may be realised by an orbifold due to Barge--Sullivan---see \cite[Theorem 13.2, p.~321]{Sul77}.

\vspace{5mm}

\textbf{ad (4)} As for the assertion on $7$-dimensional manifolds, note that every rationally elliptic Sullivan model may be realised by a compact smooth manifold in dimension not divisible by four. So it suffices to consider the product Sullivan algebra corresponding to $\s^3\times M$ with $M$ four-dimensional running over the infinitely many pairwise distinct rational types we just constructed. For degree reasons any morphism between the minimal models of these spaces has to respect the product splitting. Hence the $7$-manifolds cannot be rationally equivalent.

\vspace{5mm}

\textbf{ad (5)} Denote by $x_i$ a homogeneous basis of $\pi_\odd(M)\otimes \qq$ and by $y_i$ a homogeneous basis of $\pi_\even(M)\otimes \qq)$.
Using the equations \eqref{eqn01}, \eqref{eqn02}, \eqref{eqn03}, \eqref{eqn04} and \eqref{eqn05}
we deduce that $\sum \deg y_i\leq 7$ and $\sum \deg x_i\leq 13$. Since $M$ is simply-connected, this implies that $\deg y_i\geq 2$ and there are at most three even-degree generators $y_i$. Analogously, there are at most four odd-degree elements $x_i$.

For degree reasons and $1$-connectedness, any non-vanishing Massey product in $H^*(M)$ has to have degree $5$. This implies that $\pi_{\even}(M)\otimes \qq=\langle y_1,y_2,y_3\rangle$ projects surjectively to the generators of the cohomology algebra $H^*(M)$ under the Hurewicz map. Passing from the minimal Sullivan model $(\Lambda V,\dif)$ of $M$ to the associated pure model $(\Lambda V,\dif_\sigma)$ (which is finite-dimensional if and only if so is the model itself---see \cite[Proposition 32.4, p.~438]{FHT01}) we observe that in this case there are at least three odd degree elements $x_1,x_2,x_3$ generating three dimensional odd rational homotopy. Since $\dim M$ is odd, we even obtain that $\pi_\odd(M)\otimes\qq=\langle x_1,x_2,x_3,x_4\rangle$ of dimension $4$. Since the dimension of the associated pure model is the same as $\dim M$ (see \cite[Proposition 32.7, p.~442]{FHT01}), we may assume that $\dif(y_i)=\dif_\sigma(y_i)=0$ for $i\in\{1,2,3\}$. Thus all the $y_i$ define cohomology classes and $\deg y_i\in \{2,4\}$ by Poincar'e duality. Suppose $\deg y_i=4$ for all three of them, then, for $\dim M$ to be finite-dimensional, three of the elements $x_i$ have to lie in degrees $4k-1$ with $k\geq 2$.
The formula
\begin{align}\label{xeqn01}
7&=\sum \deg x_i - \sum (\deg y_i-1)
\end{align}
then yields a contradiction. Suppose exactly two of the $y_i$, say $y_1,y_2$, lie in degree $4$, the element $y_3$ lies in degree $2$. Then, by analogous arguments, two of the $x_i$ lie in degrees at least $5$, with at least one in degree $7$ or larger. Formula \eqref{xeqn01} again yields a contradiction---the dimension would have to be $11$ at least. (Alternatively, one might argue using Poincar\'e duality $\pi_3(M)\otimes \qq\cong H^3(M)\cong H^4(M)$.)

If $\deg y_1=\deg y_2=2$ and $\deg y_3=4$, formula \eqref{xeqn01} yields $x_1+x_2+x_3+x_4=12$ and $\deg x_i=3$ for $i\in \{1,2,3,4\}$. Thus the subalgebra generated by $[y_3]$ is infinite-dimensional; a contradiction.

If $\deg y_1=\deg y_2=\deg y_3=2$, formula \eqref{xeqn01} yields $x_1+x_2+x_3+x_4=10$; a contradiction, since $\deg x_i\geq 3$ due to $1$-connectedness.

\vspace{3mm}

Suppose next that there are exactly $2$ generators $y_1,y_2$. Then there are at least $3$ elements $x_1,x_2,x_3$, possibly also $x_4$. As above we see that $\deg y_1,\deg y_2 \in \{2,4\}$. Suppose first $\deg y_1=\deg y_2=4$. Then, again, at least two of the $x_i$ lie in degrees $4k-1$ with $k\geq 2$ and $\dim M\geq 11$; a contradiction.

Assume now $\deg y_1=2, \deg y_2=4$. Again, the finite formal dimension of $H^*(M)$ implies that $\deg x_1=4k-1$ with $k\geq 2$. Formula \eqref{xeqn01} applied to $\deg x_1\geq 7, \deg x_2\geq 3, \deg x_3 \geq 3$ yields a contradiction.

Suppose $\deg y_1=\deg y_2=2$. Formula \eqref{xeqn01} implies that $\dim \pi_{\odd}(M)\otimes \qq=3$ and $\deg x_1=\deg x_2=\deg x_3=3$. If, without restriction, $\dif x_1=0$, then $M$ rationally splits as a product of $\s^3\simeq_\qq (\Lambda \langle x_1\rangle,0)$ and a four-dimensional rationally elliptic space with rational cohomology generated by $[\langle y_1,y_2\rangle]$. The classification of these algebras from above
shows that $M$ has the rational homotopy type of $\s^3\times (\cc\pp^2\# \cc\pp^2)$ or of $\s^3\times \s^2\times \s^2$ in this case.

If $\dif x_1\neq 0$, i.e.~if $\dif|_{\langle x_1,x_2,x_3\rangle}$ is injective, we obtain that $H^4(M)=0$ and this case uniquely corresponds to the rational bundle of $\s^3$ over $\s^2\times \s^2$ where $\s^3$ is attached to the volume form up to a non-trivial multiple.

Assume now there is only one even-degree generator $y_1$. The finite-dimensionality of the cohomology algebra implies that $\dif x_1=y_1^k$ for some $k\geq 2$. Thus $M$ rationally splits as a product of a rational $\cc\pp^{k-1}$ or a rational $\hh\pp^{k-1}$ and odd degree sphere factors. Formula \eqref{xeqn01} together with $1$-connectedness then implies the possibilities from the assertion.

If the rational homotopy is concentrated in odd degrees only, $1$-connectedness directly shows that $M$ is a rational sphere.

It is trivial to see from homotopy groups combined with obvious arguments using cohomology that the given spaces are not rationally equivalent. In the case of $\s^3\times \s^2\times\s^2$ and $\s^3\times \cc\pp^2\#\cc\pp^2$ one uses the argument from (4) again.
\end{prf}

\begin{proof}[\textsc{Proof of Theorem \ref{theo03}}]
The first assertion on the non-existence of metrics of holonomy $\Spin(7)$ can be derived as follows:
Relation \eqref{PUReqn08} yields $b_4\geq b_4^+\geq b_2-b_3+25$ and
$b_3+b_4\geq b_2+25$. Let $(\Lambda V,\dif )$ be the
minimal model of $M$. Since $M$ is simply-connected, we obtain that
\begin{align*}
\dim \pi_2(M)\otimes \qq&=\dim V^2=b_2\\
\dim \pi_3(M)\otimes \qq&=\dim V^3\geq b_3\\
\dim \pi_4(M)\otimes \qq& =\dim V^4\geq b_4-\dim
\Sym_2(V^2)=b_4-\frac{b_2(b_2+1)}{2}
\end{align*}
In particular, we see that
\begin{align*}
\dim \pi_*(M)\otimes \qq&\geq
b_2+b_3+\bigg(b_4-\frac{b_2(b_2+1)}{2}\bigg)\\&=
b_3+b_4+\frac{b_2(1-b_2)}{2}\\&\geq
25+b_2+\frac{b_2(1-b_2)}{2}\\&=25+\frac{b_2(3-b_2)}{2}
\end{align*}
Since we have that $\dim V^2=b_2$, equation \eqref{eqn02} yields,
in particular, that \linebreak[4]$2\cdot b_2\leq \dim M=8$ and that $b_2\leq 4$.
Thus we derive that $\dim \pi_*(M)\otimes \qq\geq 23$, which
contradicts \eqref{eqn03}.

\vspace{3mm}

The assertion on holonomy contained in $\SU(4)$ follows from Inequality \eqref{PUReqn10}, which implies that $b_3+b_4\geq 50$. Indeed, arguing as before using \eqref{eqn02}, the Betti numbers $b_2$ and $b_3$ are restricted from above by $\dim \pi_2(M)\otimes \qq\leq 4$ and $\dim \pi_3(M)\otimes \qq\leq 5$ respectively. The fourth Betti number is restricted from above by $\dim \pi_4(M)\otimes \qq +\frac{b_2(b_2+1)}{2}\leq 2+10=12$.
It follows that $b_3+b_4\leq 5+12=17$ and the manifold $M$ cannot carry holonomy $\SU(4)$.

From \cite[Theorem 10.6.1, p.~259]{Joy00} we derive that $\Hol(M)=\Sp(2)$ if and only if $\hat A(M)=3$. Due to \cite[Corollary 3.6.3, p.~67]{Joy00} the manifold $M$ is spin. Thus a homogeneous structure implies in particular the existence of a smooth $\s^1$-action on $M$. Due to \cite{AH70} the $\hat A$-genus has to vanish in this case; a contradiction.

\vspace{3mm}

From Proposition \ref{lemma02} we directly see that the cohomological properties of a manifold with $\G_2$-holonomy can only be satisfied by the real types of $\s^4\times \s^3$, $\cc\pp^2\times \s^3$, $\s^3\times \cc\pp^2\# \cc\pp^2$. Indeed, the only types with $b_3\neq 0$ are these and additionally $\s^2\times \s^2\times \s^3\simeq_\qq\s^3\times \cc\pp^2\#\overline{\cc\pp^2}$.

This latter case can be excluded by the structure of the cohomology algebra. We use the description as $\s^2\times \s^2\times \s^3$ and we denote the two generators of the second rational cohomology corresponding to the respective $\s^2$-factors by $a$ and $b$. We compute $a^2\cdot \omega=0$ contradicting \eqref{PUReqn09}.

It is obvious that any of these real types is formal---recall that formality does not depend of the extension field of $\qq$. Due to Berger and \cite[Theorem 10.5.7, p.~256]{Joy00} an irreducible manifold with holonomy contained in $\Spin(7)$ has either holonomy $\Spin(7)$, $\G_2$, $\SU(4)$, $\Sp(2)$ or is a symmetric space. (Note that the holonomy group is necessarily connected.) In either case (using the formality of K\"ahler manifolds from \cite{DGMS75}) it is formal once we assume it to be rationally elliptic. This finishes the proof.
\end{proof}


\section{Proof of Theorem \ref{theo01}}\label{sec03}

Due to \cite[Lemma 3.3, p.~408]{Tot02} we may assume throughout this article (up to diffeomorphism) that a biquotient $\biq{G}{H}$ satisfies that $G$ is simply-connected, $H$ is connected and $H$ does not act transitively on any simple factor of $G$.

Denote by $d(G)$ the largest degree of a non-vanishing rational homotopy group of the compact Lie group $G$. This degree is given as in Table \ref{table01}---see \cite[Table 4.2, p.~410]{Tot02}.
\begin{table}[h]
\centering \caption{degrees of rational homotopy groups}
\label{table01}
\begin{tabular}{c@{\hspace{8mm}}|  @{\hspace{8mm}}c@{\hspace{8mm}}|@{\hspace{8mm}}c@{\hspace{8mm}}}
$G$ & $i$ with $\dim \pi_i(G)\otimes \qq=1$ &$d(G)$ \\
\hline
$\SU(n)$,&$3,5,7,\dots, 2n-1$ & $2n-1$\\
$\SO(n)$, $n\geq 4$ even& $3,7,11,\dots,2n-5,n-1$ & $2n-5$\\
$\SO(n)$, $n$ odd& $3,7,11,\dots,2n-3$& $2n-3$\\
$\SO(2)$ &$1$& $1$\\
$\Sp(n)$, &$3,7,11,\dots,4n-1$& $4n-1$  \\
$\G_2$ &$3,11$& $11$\\
$\F_4$ &$3,11,15,23$& $23$\\
$\E_6$ &$3,9,11,15,17,23$& $23$\\
$\E_7$ &$3,11,15,19,23,35$& $35$\\
$\E_8$ &$3,15,23,27,35,39,47,59$& $59$
\end{tabular}
\end{table}
We cite \cite[Lemma 4.3, p.~410]{Tot02}.
\begin{lemma}\label{lemma01}
Let $H\to G$ be a non-trivial homomorphism of simply-connected simple Lie groups. Then the maximal degrees of non-vanishing rational homotopy groups $d(H)$, $d(G)$ satisfy $d(H)\leq d(G)$ with equality if and only if either the homomorphism is bijective or $G/H$ is one of the following homogeneous spaces.
\begin{align*}
&\s^{2n-1} \textrm{ for } n\geq 4, \quad \SU(2n)/\Sp(n) \textrm{ for } n\geq 2, \quad \\&\s^7\times \s^7=\Spin(8)/\G_2, \quad\E_6/\F_4
\end{align*}
The case of $\s^{2n-1}$ is realised by $\Spin(2n)/\Spin(2n-1)$ and by $\Spin(7)/\G_2=\s^7$.
\end{lemma}

Recall that we assumed $G$ to be simply-connected and $H$ to be connected.
Up to finite covering we may split $H=H_1\times \ldots \times H_k$ with simple simply-connected factors and $\s^1$-factors. We also split $G=G_1\times \ldots G_l$ into simple factors. Up to diffeomorphism we may assume that $H$ does not act transitively on any factor of $G$---see \cite[Corollary 4.6]{Tot02}.

Note that Proposition \ref{lemma02} implies that the homotopy Euler characteristic
\begin{align*}
\chi_\pi(M)=\dim \pi_\odd(M)-\dim \pi_\even(M)=1
\end{align*}
for a simply-connected seven-dimensional rationally elliptic space $M$. Since $\chi_\pi(\biq{G}{H})=\rk G - \rk H$ we derive that
\begin{align}\label{xeqn02}
\rk G= \rk H+1
\end{align}

The proof of Theorem \ref{theo01} will basically proceed in two cases. In the first case we assume that for all $i,j$ we have that the orbit inclusion $H\hto{} G$ satisfies that whenever there is an injective projection $H_i\hto{} G_j$ we have $d(H_i)<d(G_j)$. In the second case we assume that there exists an $H_i$ mapping injectively into some $G_j$ with $d(H_i)=d(G_j)$.

\subsection{Case 1}
 $d(H_i)<d(G_j)$. Depending on the rational type of the biquotient---using the top degree of its rational homotopy groups---we determine all potential Lie groups (up to isomorphisms and finite coverings) which may be factors $G_j$ in Table \ref{table02}.
\begin{table}[h]
\centering \caption{potential factors $G_j$}
\label{table02}
\begin{tabular}{c@{\hspace{4mm}}| @{\hspace{4mm}}c@{\hspace{4mm}} |@{\hspace{4mm}}c@{\hspace{4mm}}}
real type & $\deg \pi_*(M)$ &$G_j$ \\
\hline
$\s^7$ & $7$& $\SU(4), \SU(3), \SU(2), \SO(5)$\\
$\s^4\times \s^3$ & $3,4,7$ & $\SU(4), \SU(3), \SU(2), \SO(5)$\\
$\s^2 \times \s^5$ & $2,3,5$ & $\SU(3), \SU(2)$\\
$\cc\pp^2 \times \s^3$ &$2,3,5$ &$\SU(3), \SU(2)$\\
remaining & $2,2,3,3,3$ & $\Sp(1)$
\end{tabular}
\end{table}

In Table \ref{table03} we determine all potential pairs $(G,H)$.

 This table arises as follows:
\begin{itemize}
\item
First, we determine the maximal number $l$ of factors of $G$. For this we consider the orbit map $H\to{} G_j$ for the respective factors of $G$. If the induced map on rational homotopy groups $\pi_*(H)\otimes \qq\to \pi_*(G_j)\otimes \qq$ is surjective, then $H$ acts transitively on $G_j$ and, up to diffeomorphism, we may cancel this factor---see \cite[Corollary 4.6]{Tot02} and \cite[Lemma 3.2]{Tot02}. Thus, we may assume that each factor $G_j$ non-trivially contributes rational homotopy groups to $\biq{G}{H}$. By the construction of the Sullivan model of this biquotient, homotopy groups contributed by $G$ are necessarily of odd degree. Thus the respective rational homotopy types provide effective upper bounds on the number of factors of $G$. More precisely, in the case of $\s^7$ there is only one factor, in the case of $\s^4\times \s^3$ there are at most two factors, for $\s^2\times \s^5$ also at most two factors, as well as for $\cc\pp^2\times \s^3$, in the remaining cases there are up to three factors.
\item
Second, we use Property \eqref{xeqn02} and the fact that $7=\dim \biq{G}{H}=\dim G-\dim H$ in order to determine all potential groups $G$. Note that only one factor may have a non-trivial top degree rational homotopy group; so there has to be exactly one such factor.
\item We do have to compute rational homotopy groups and Betti numbers in some cases to prove that certain pairs $(G,H)$ cannot realise a biquotient. An example of this is $(\SU(3)\times \SU(3),\Sp(1)\times \Sp(1)\times \Sp(1))$. This cannot realise a biquotient, since there will be a generator of the rational cohomology algebra in degree $4$ which will generate an algebra isomorphic to $H^*(\hh\pp^\infty)$.
\item Finally, this yields Table \ref{table03}. Note that the pairs
\begin{align*}
&(\SU(4)\times \SU(3), \SU(3)\times \SU(3))\\
&(\SU(4)\times \Sp(1),\SU(3)\times \SU(2))\\
&(\Sp(2)\times \SU(3), \Sp(1)\times \SU(3))\\
&(\SU(3)\times \Sp(1),\s^1\times \Sp(1))
\end{align*}
a priori might arise from orbit inclusions which do not map the second factor of $H$ onto $G_2$.
\end{itemize}

\begin{table}[h]
\centering \caption{potential pairs $(G,H)$, Case 1}
\label{table03}
\begin{tabular}{c@{\hspace{4mm}}| @{\hspace{4mm}}c@{\hspace{4mm}}}
real type & $(G,H)$ \\
\hline
$\s^7$ & $(\SU(4),\SU(3)$, $(\Sp(2),\Sp(1))$, $(\SU(3),\s^1)$ \\
$\s^4\times \s^3$ &  $(\SU(4), \SU(3))$, $(\Sp(2),\Sp(1)$ , $(\SU(3),\s^1)$,\\& $(\SU(4)\times \SU(3), \SU(3)\times \SU(3))$, \\&$(\SU(4)\times \Sp(1),\SU(3)\times \SU(2))$, \\&$(\Sp(2)\times \SU(3), \Sp(1)\times \SU(3))$\\
$\s^2 \times \s^5$ &  $(\SU(3),\s^1)$, $(\SU(3)\times \Sp(1),\s^1\times \Sp(1))$\\
$\cc\pp^2 \times \s^3$ &$(\SU(3),\s^1)$, $(\SU(3)\times \Sp(1),\s^1\times \Sp(1))$\\
remaining &  $(\Sp(1)\times \Sp(1)\times \Sp(1),\s^1\times \s^1)$
\end{tabular}
\end{table}

Let us now provide the following lemma to reduce the list of potential biquotients.

\begin{lemma}
\begin{itemize}
\item
We have $M=\biq{\SU(3)}{\s^1}\simeq_\qq \s^2\times\s^5$ (irrespective of the action of the $\s^1$).
\item
A biquotient $\biq{\Sp(2)}{\Sp(1)}$ is either diffeomorphic to the $7$-sphere or the Gromoll--Meyer sphere.
\item
A biquotient $\biq{\SU(4)}{\SU(3)}$ is diffeomorphic to the $7$-sphere.
\end{itemize}
\end{lemma}
\begin{prf}
The first assertion can be derived as follows.
Suppose the $\s^1$ acts via
\begin{align*}
(\diag(z^{a_1},z^{b_1},z^{-a_1-b_1}),\diag(z^{a_2},z^{b_2},z^{-a_2-b_2}))
\end{align*}
We build the Sullivan model of this biquotient as
\begin{align*}
(\Lambda \langle u,x,y \rangle, \dif)
\end{align*}
with $\deg u=2, \deg x=3, \deg y=5$ and $\dif u=0$,
\begin{align*}
\dif x=-a_1^2 + 2 a_1 a_2 - a_2^2 - a_1 b_1 + a_2 b_1 - b_1^2 + a_1 b_2 - a_2 b_2 +
 2 b_1 b_2 - b_2^2
\end{align*}
It is obvious that $M$ has not the rational type of $\s^2\times \s^5$ if and only if $\dif x=0$. This yields
\begin{align*}
a_1=\frac{1}{2} \bigg(2 a_2 - b_1 + b_2 \pm \sqrt{-3(b_1-b_2)^2}\bigg)
\end{align*}
Thus we deduce $b_1=b_2$ and $a_1=a_2$ in this case. Yet, the $\s^1$-action given by these rotation numbers is not free.

\vspace{3mm}
For the second assertion we use the classification of homogeneous $7$-manifolds (see \cite{Kla88}) in order to see that there is no such space which has a non-trivial third rational homotopy group, but which is rationally $2$-connected. In other words, this implies that the orbit inclusion $\Sp(1)\hto{}\Sp(2)$ of the biquotient must be an isomorphism on the rationalised third homotopy group. Consequently, $\pi_*(\biq{\Sp(2)}{\Sp(1)})\otimes \qq$ is concentrated in degree $7$ and a rational $7$-sphere. Due to the classification in \cite{KZ04} we derive that this biquotient is either the sphere or the Gromoll--Meyer sphere.

\vspace{3mm}

The third assertion can be deduced easily, since the only inclusion (up to conjugation) of $\SU(3)$ into $\SU(4)$ is the standard one. So the orbit inclusion is rationally $6$-connected. Again by \cite{KZ04} we see that this biquotient is diffeomorphic to the $7$-sphere.
\end{prf}

For rank reasons there is no inclusion $\SU(3)\times \SU(3)\hto{}\SU(4)$ and thus the projection to $G_2=\SU(3)$ has to be non-trivial. Since this is a simple Lie group, we obtain the following:
The induced orbit map for the pair $(\SU(4)_\qq\times \SU(3)_\qq, \SU(3)_\qq\times \SU(3)_\qq)$  (where the subscript denotes rationalisation) has to be surjective after projection onto the factor $G_2$, i.e.~$\pi_*(\SU(3)\times \SU(3))\otimes \qq\to \pi_*(\SU(4)\times \SU(3))\otimes \qq \to \pi_*(\SU(3))\otimes \qq$ is surjective and, up to diffeomorphism the pair $(\SU(4)\times \SU(3), \SU(3)\times \SU(3))$ corresponds to $(\SU(4),\SU(3))$.

A similar argument applies to the pair $(\SU(4)\times \Sp(1),\SU(3)\times \SU(2))$. Indeed, there is no inclusion of $\SU(3)\times \SU(2)$ into $\SU(4)$ so that the action is transitive on $G_2=\Sp(1)$ and both $\Sp(1)\cong \SU(2)$ factors can be cancelled.

Since $\SU(3)$ cannot be realised as a subgroup of $\Sp(2)$ (not even up to finite covering), an analogous argument also applies to $(\Sp(2)\times \SU(3), \Sp(1)\times \SU(3))$.

\begin{lemma}\label{lemma04}
A biquotient $\biq{\SU(3)\times \Sp(1)}{\s^1\times \Sp(1)}$ has either the rational type of $\s^2\times \s^5$ or is diffeomorphic to the product $\cc\pp^2\times \s^3$.
\end{lemma}
\begin{prf}
If the orbit inclusion does not map $H_2=\Sp(1)$ surjectively to $G_2=\Sp(1)$ (up to finite covering and after projecting to $G_2$), then we have to differentiate between the two cases: Either $H_1=\s^1$ maps non-trivially into $G_2=\Sp(1)$ or not. In the first case, the model of this biquotient yields that the generator of the cohomology algebra of the biquotient corresponding to this $\s^1$ has vanishing square. In the second case, the biquotient is just $\cc\pp^2\times \s^3$. Indeed, in this case we have an orbit inclusion $\s^1\times \Sp(1)\hto{} \SU(3)$. Up to conjugation and due to the classification of maximal rank subgroups this is the standard inclusion $\mathbf{S}(\U(2)\times \U(1))\hto{}\SU(3)$. Thus the action of $\s^1\times \Sp(1)$ is trivial on $G_2=\Sp(1)$. Thus the biquotient splits as a product of $\Sp(1)$ and an equal rank biquotient of $\SU(4)$ of the rational type of $\cc\pp^2$. According to the classification of biquotients with singly generated ratonal cohomology algebra (see \cite{KZ04}) we obtain that $\biq{H_1}{\s^1\times \Sp(1)}=(\U(3)/\U(1)\times \U(2))\times \Sp(1)=\cc\pp^2\times \Sp(1)$.
\end{prf}

\begin{lemma}
A biquotient $\biq{\Sp(2)\times \Sp(1)}{\Sp(1)\times \Sp(1)}$ is either diffeomorphic to the sphere, the Gromoll--Meyer sphere or
$\s^4\times \s^3$.
\end{lemma}
\begin{prf}
Consider the orbit inclusion $\Sp(1)\times \Sp(1) \hto{} \Sp(2)\times \Sp(1)$. Either the induced map on $\pi_3(\cdot)\otimes \qq$ is surjective on $G_2=\Sp(1)$ and the biquotient is diffeomorphic to $\biq{\Sp(2)}{\Sp(1)}$ or the orbit includes completely into $G_1=\Sp(2)$. In this case the biquotient splits as $\biq{Sp(2)}{\Sp(1)\times \Sp(1)}$. In either case we use the classification in \cite{KZ04} to see that we either obtain the sphere, the Gromoll--Meyer sphere (first case) or $\s^4\times \s^3$ (second case) up to diffeomorphism.
\end{prf}

These arguments combine to see that Table \ref{table04} provides a complete list of biquotients. It is trivial to see that every given real homotopy type can be realised by one of these biquotients.
For the case of the non-trivial rational bundle we use
\begin{cor}
A simply-connected rationally elliptic space is formal unless it has the real type of $\s^3\tilde\times (\s^2\times\s^2)$. This type can be realised by the homogeneous space $\Sp(1)\times \Sp(1)\times \Sp(1)/\s^1\times\s^1$ with the inclusion $(a,b)\mapsto (a,b,a\cdot b)$ already.
\end{cor}
\begin{prf}
This follows from Proposition \ref{lemma02} and a simple computation of the rational homotopy type of the given homogeneous space.
\end{prf}

\begin{table}[h]
\centering \caption{biquotients $\biq{G}{H}$, Case 1}
\label{table04}
\begin{tabular}{c@{\hspace{4mm}}| @{\hspace{4mm}}c@{\hspace{4mm}}}
real type & $\biq{G}{H}$ \\
\hline
$\s^7$ & $\biq{\SU(4)}{\SU(3)}$, $\biq{\Sp(2)}{\Sp(1)}$ \\
$\s^4\times \s^3$ & $(\Sp(2)/\Sp(1)) \times \Sp(1)$\\
$\s^2 \times \s^5$ &  $\biq{\SU(3)}{\s^1}$, $\biq{\SU(3)\times \Sp(1)}{\s^1\times \Sp(1)}$\\
$\cc\pp^2 \times \s^3$ &$\biq{\SU(3)\times \Sp(1)}{\s^1\times \Sp(1)}$\\
remaining &  $\biq{\Sp(1)\times \Sp(1)\times \Sp(1)}{\s^1\times \s^1}$
\end{tabular}
\end{table}

\vspace{5mm}

\subsection{Case 2} Suppose now that in the orbit inclusion there is a factor $H_i$ mapping to some $G_j$ with $d(H_i)=d(G_j)$. Without restriction, we assume $i=j=1$.

We say that a homotopy group of $G$  in degree $i$ survives to $\biq{G}{H}$ if the map $\pi_i(G)\to \pi_i(\biq{G}{H})$ induced by the projection is injective. (Analogously for rational homotopy groups.)
\begin{lemma}\label{lemma03}
Let $M^7=\biq{G}{H}$ be a simply-connected biquotient. Then either the homotopy group of $G$ of largest degree or the one of second largest degree survives to $\biq{G}{H}$.
\end{lemma}
\begin{prf}
We show that if the largest degree does not survive, the second largest one does. In this case $G_1$ and $H_1$ are as in Lemma \ref{lemma01}.

If $G_1/H_1=\s^{2n-1}$ with $n\geq 4$, then, in order to kill the largest remaining degree, there must be a factor $H_2$ which is mapped under the orbit inclusion into $G_1$ and such that $\pi_{2n-1}(H_2)\otimes \qq$ maps surjectively to $\pi_{2n-1}(G_1)\otimes \qq$. For rank reasons---$H_1\times H_2$ includes into $G_1$ then---such an $H_2$ does not exist.

In the case $G_1/H_1=\SU(2n)/\Sp(n)$ with $n\geq 2$, we would need an $H_2$ with $\rk H_2\leq n-1$ including into $\SU(2n)$ and killing degree $4n-3$ homotopy. For a classical Lie group it is easy to see that this is not possible---the case of largest degree rational homotopy with smallest rank being realised by $\Sp(n)$. For exceptional Lie groups we use Table \ref{table01}. The only Lie group which might kill an element of degree $4n-3$ is $E_6$ with degrees $9$ and $17$. However, $4n-3=9$ implies $n=3$, $4n-3=17$ implies $n=5$ and $\E_6\not\subseteq \SU(6)$, $\E_6\times \Sp(5)\not\subseteq \SU(10)$.

If $G_1/H_1=\Spin(7)/G_2=\s^7$, we cannot find a group of rank $1$ killing an element of degree $7$.

If $G_1/H_1=\Spin(8)/\G_2=\s^7\times \s^7$, the only group $H_2$ of rank at most two potentially killing an element in degree $7$ is $\Sp(2)$. However, $\G_2\times \Sp(2)$ is not a subgroup of $\Spin(8)$ by the classification of maximal rank subgroups due to \cite{BD49}.

In the case $G_1/H_1=\E_6/\F_4$ it is impossible to kill an element of degree $17$ by a group of rank at most $2$.
\end{prf}

We use this lemma to produce Table \ref{table05} of potential factors of $G$ in Case 2 (up to isomorphism and finite covering). That is, we identify all simple groups which have either top degree rational homotopy smaller than the top degree rational homotopy of the respective case or which appear in Lemma \ref{lemma01} and have the property that their second largest rational homotopy group lies in a degree not larger than the largest one permitted by the respective rational homotopy type.

Since we are in Case 2, we may assume that at least one factor $G_1$ is one of
\begin{align}\label{xeqn03}
\SU(8), \SU(6),\SU(4), \SO(8), \SO(7)
\end{align}
due to Lemma \ref{lemma01}.

\begin{table}[h]
\centering \caption{potential factors $G_j$, Case 2}
\label{table05}
\begin{tabular}{c@{\hspace{4mm}}| @{\hspace{4mm}}c@{\hspace{4mm}} |@{\hspace{4mm}}c@{\hspace{4mm}}}
real type & $\deg \pi_*(M)$ &$G_j$ \\
\hline
$\s^7$ & $7$&  $\SU(4), \SU(3), \SU(2), \SO(8),$\\&&$ \SO(7), \SO(5)$\\
$\s^4\times \s^3$ & $3,4,7$ & $\SU(4), \SU(3), \SU(2), \SO(8),$\\&&$ \SO(7), \SO(5)$\\
$\s^2 \times \s^5$ & $2,3,5$ & $\SU(4),\SU(3), \SU(2)$\\
$\cc\pp^2 \times \s^3$ &$2,3,5$ &$\SU(4),\SU(3),\SU(2)$\\
remaining & $2,2,3,3,3$ & $\SU(3), \SU(2)$
\end{tabular}
\end{table}

We have the same maximal number of factors $G_j$ in the respective cases as above, we use $\rk G=\rk H+1$ and $\dim G=7+\dim H$ in order to determine all potential pairs $(G,H)$ with $G_1$ from the list in \eqref{xeqn03} which may arise in Case 2. Combining these general arguments with concrete computations of rational homotopy groups and Betti numbers using concrete potential inclusions of Lie groups we see that $\biq{G}{H}$ is as in Table \ref{table06}.
\begin{table}[h]
\centering \caption{biquotients $\biq{G}{H}$, Case 2}
\label{table06}
\begin{tabular}{c@{\hspace{4mm}}| @{\hspace{4mm}}c@{\hspace{4mm}}}
real type & $\biq{G}{H}$ \\
\hline
$\s^7$ & $\Spin(7)/\G_2, \SO(8)/\SO(7)$ \\
$\s^4\times \s^3$ & \\
$\s^2 \times \s^5$ &  $\biq{\SU(4)\times \Sp(1)}{\Sp(2)\times \s^1}$\\
$\cc\pp^2 \times \s^3$ &\\
remaining &
\end{tabular}
\end{table}
We used
\begin{lemma}
We have $\biq{\SU(4)\times \Sp(1)}{\Sp(2)\times \s^1}\simeq_{\qq}\s^2 \times \s^5$.
\end{lemma}
\begin{prf}
Since $\Sp(2)\times \s^1$ is not a maximal rank subgroup of $\SU(4)$, the orbit inclusion maps $\s^1$ non-trivially to $\Sp(1)$. That is, identifying the biquotient up to diffeomorphism with $H\backslash G\times G/\Delta G$ (cf.~\cite{Sin93}) this map is identified with $H\to G\times G/\Delta G$. Thus the projection of the inclusion $\s^1\hto{} (\SU(4)\times \Sp(1))^2$ onto both $\Sp(1)$-factors cannot be identical.

Moreover, there is no non-trivial map from $\Sp(2)$ to $\Sp(1)$. Building the model of the biquotient we thus see that the square of the generator of the rational cohomology algebra corresponding to the $\s^1$-factor vanishes. Thus the rational homotopy type of the biquotient is $\s^2\times \s^5$.
\end{prf}

Now both Tables \ref{table04} and \ref{table06} combine to yield Table \ref{table07} and prove Theorem \ref{theo02}.

\vspace{5mm}

Let us now prove Theorem \ref{theo01}. We need the following lemma, which follows the arguments in \cite[Example 8.4, p.~130]{Cav05}.
\begin{lemma}\label{lemma05}
Suppose $M^7=M^4\times \s^3$. Then $M^7$ cannot carry a metric of holonomy $\G_2$.
\end{lemma}
\begin{prf}
Since $M^7$  is simply-connected, so is $M^4$ and $H^3(M^7)=H^3(\s^3)$ generated by a pullback of the volume form $\omega$ of $\s^3$. Thus the bilinear form $(a,b)\mapsto \int_{M^7} a\cdot b \cdot \omega$ (with $a,b,\in H^2(M^7)$) is just the intersection form on $M^4$. If $M^7$ carries a metric of holonomy $\G_2$, the intersection form of $M^4$ is (negative) definite. Since $M^7$ is spin, so is $M^4$. Due to Donaldson (cf.~\cite[Theorem 1.3.1]{DK90}) the manifold $M^4$ can only be spin, if $b_2(M^4)=0$. In this case the Hirzebruch signature formula implies that $p_1(M^4)=3\sigma(M^4)=0$ where $\sigma(M^4)$ denotes the signature of $M^4$. This implies that both $M^4$ and then also $M^7$ have vanishing first Pontryagin class; a contradiction.
\end{prf}

\begin{proof}[\textsc{Proof of Theorem \ref{theo01}}]
From Table \ref{table07} we see that a simply-connected biquotient has $b_3\neq 0$ if an only if it is
$\s^4\times \s^3$, $\biq{\SU(3)\times \Sp(1)}{\s^1\times \Sp(1)}$ or $\biq{\Sp(1)\times \Sp(1)\times \Sp(1)}{\s^1\times \s^1}$. Due to Lemma \ref{lemma04} we have that a biquotient $\biq{\SU(3)\times \Sp(1)}{\s^1\times \Sp(1)}$ has either the rational type of $\s^2\times \s^5$ or is diffeomorphic to the product $\cc\pp^2\times \s^3$. In total, this means that unless $M^7$ splits as a genuine product $M^4\times \s^3$ up to diffeomorphism, $M$ is  $\biq{\Sp(1)\times \Sp(1)\times \Sp(1)}{\s^1\times \s^1}$. The product manifold cannot carry $\G_2$-holonomy due to Lemma \ref{lemma05}.

For the case of homogeneous spaces we either observe that also the remaining biquotient case then can only be realized by a genuine product or we use the classification in \cite{Kla88} together with the computations of the Betti numbers there. They show exactly that a seven dimensional homogeneous space can only have $b_3\neq 0$ if it splits as a product with an $\s^3$-factor. Alternatively, we might use that a homogeneous space of the form $G/T$ is stably parallelisable; hence its first Pontryagin class vanishes.
\end{proof}


\section{On properties of positive quaternion K\"ahler manifolds}\label{sec04}

Let us now deal with closed oriented manifolds of holonomy $\Sp(n)\Sp(1)$, quaternion K\"ahler manifolds. They are Einstein and thus their scalar curvature is either negative, positive or vanishes. A complete quaternion K\"ahler manifold of positive scalar curvature is called a \emph{positive quaternion K\"ahler manifold}. Due to LeBrun--Salamon they are conjectured to be symmetric and various classification results do confirm this in special cases (see \cite{Sal82}, \cite{Sal99}, \cite{Fan08}, \cite{AD10}, \cite{Ama12},  etc.) In the symmetric case (and more generally in the homogeneous case) there is a classification yielding one such space for each of the classical and the exceptional complex Lie algebras. The symmetric positive quaternion K\"ahler manifolds, the so-called \emph{Wolf spaces}, corresponding to the exceptional Lie algebras are called the \emph{exceptional Wolf spaces} (like $\G_2/\SO(4)$, \dots).

A positive quaternion K\"ahler manifold $M$ is simply-connected. We denote by $b_i=\dim H^i(M;\qq)$ the Betti numbers of $M$, by $c_i=\dim \pi_*(M)\otimes \qq$ the homotopy Betti numbers and by $(\Lambda V,\dif)$ a minimal model of $M$.

We may assume that $b_2=0$, since otherwise $M$ is homothetic to the symmetric spaces $\Gr_2(\cc^{n+2})$---see \cite[Theorem 5.5, p.~103]{Sal99}. Odd Betti numbers of $M$ vanish (see \cite[Theorem 6.6, p.~163]{Sal82}). In particular, we may assume $M$ to be rationally $3$-connected.
From \cite[Theorem 5.4, p.~403]{Sal93} we cite the following relations for the Betti numbers of a positive quaternion K\"ahler manifold of dimensions $12$ respectively $16$.
In the first case we have that $b_6=0$, in the second one it follows that $-1 +3 b_4-b_6=2b_8$. A positive quaternion K\"ahler manifold has the Hard-Lefschetz property with respect to a $4$-form, the Kraines form (as follows easily from the classical Hard-Lefschetz property of its twistor space, which is a K\"ahler manifold).

Positive quaternion K\"ahler manifolds exist only in dimensions divisible by four and are classified up to dimension eight (see \cite{Sal82}, \cite{AD10}).

\vspace{5mm}

Due to work of Galicki, Salamon and the author (see \cite{Sal99}, \cite{Ama11a}) it was obtained that a positive quaternion K\"ahler manifold of dimension $4n$ with $n\in\{1,3,4,5,6\}$ is homothetic to a quaternionic projective space, if its fourth Betti number equals one. (Note that the eight-dimensional exceptional Wolf space $\G_2/\SO(4)$ has the rational homotopy type of $\hh\pp^2$.)

We vary this theorem for positive quaternion K\"ahler manifolds over an underlying biquotient. (Note that both structures need not be related.)
\begin{theo}\label{theo04}
Let $M$ be a biquotient which also bears the structure of a positive quaternion K\"ahler manifold. Suppose $M$ has the rational homotopy type of a compact rank one symmetric space.

Then $M$ is homothetic to a quaternionic projective space, the exceptional Wolf space $\G_2/\SO(4)$ or the Wolf space $\Gr_2(\cc^{4})$.
\end{theo}
\begin{prf}
A positive quaternion K\"ahler manifold $M$ is simply-connected. The simply-connected compact rank one symmetric spaces are $\s^n$, $\cc\pp^n$, $\hh\pp^n$, $\operatorname{Ca}\pp^2$. As stated above, positive quaternion K\"ahler manifolds satisfy Hard-Lefschetz with respect to the Kraines form in degree $4$. Moreover, if $b_2(M)\neq 0$, then $M\cong \Gr_2(\cc^{n+2})$. It is easy to check that $\Gr_2(\cc^{n+2})\not\simeq_\qq \cc\pp^{2n}$ in dimensions larger than $4$ and that $\Gr_2(\cc^{4})\simeq_\qq \cc\pp^{2}$. Thus for cohomological reasons it follows that we can focus on the case when $M\simeq_\qq \hh\pp^n$.

We may assume that $\dim M>8$ due to the classification of positive quaternion K\"ahler manifolds in low degrees. In other words, it remains to show that for $\dim M\geq 12$, $M$ is homothetic to $\hh\pp^n$.

We apply the classification of simply-connected biquotients in \cite{KZ04}. Indeed, we see that any such biquotient which has the rational homotopy type of a quaternionic projective space is either diffeomorphic to a compact rank one symmetric space or to a space from \cite[Table I, p.~150]{KZ04}.

If $M$ is diffeomorphic to $\hh\pp^n$, it is homothetic to $\hh\pp^n$, since $\pi_2(M)=0$ determines the homothety type---see \cite[Theorem 5.5 (i), p.~103]{Sal99}.

The only biquotient $M$ from the table which satisfies $\dim M\geq 12$ and which has the rational type of a quaternionic projective space is the biquotient $\Delta \SU(2)\backslash \SO(4n+1) /\SO(4n-1)\simeq_\qq\hh\pp^{2n-1}$. 

Suppose $M$ is this biquotient. The computation on \cite[p.~158]{KZ04} shows that $H^2(M;\zz)=0$. Due to \cite[Theorem 6.3, p.~160]{Sal82} this yields that $M$ would have to be homothetic to the quaternionic projective space; a contradiction.
\end{prf}

\begin{theo}\label{LOWtheo07}
A rationally elliptic $16$-dimensional positive quaternion K\"ahler manifold is either homothetic to $\hh\pp^4$ or $\Gr_2(\cc^6)$ or has the Betti numbers and homotopy Betti numbers of $\widetilde \Gr_4(\rr^8)$.
\end{theo}
\begin{prf}
We may assume that $M$ is rationally $3$-connected unless $M\cong \Gr_2(\cc^6)$.
Since $M$ is an $F_0$-space and formal, we cite
\begin{align}\label{eqn09}
\dim V^\even=\dim V^\odd\leq \cat(\Lambda V,\dif)=c(\Lambda V,\dif)=4
\end{align}
from \cite[Theorem 32.6 (iv), p.~441]{FHT01} and \cite[Example 4, p.~388]{FHT01} using Hard-Lefschetz. (Here $\cat$ denotes Lusternik--Schnirelmann category and $c$ the cup-length.)

We combine the equation
\begin{align*}
-1+3b_4-b_6&=2b_8
\end{align*}
with the Hard-Lefschetz property to obtain that
\begin{align*}
(b_4,b_6,b_8)\in \{&(1,0,1),
(2,1,2), (3,0,4), (3,2,3), (4,1,5), (4,3,4), (5,0,7),\\&
(5,2,6), (5,4,5), \dots \}
\end{align*}
(The remaining triples satisfy $b_4\geq 6$.)
It follows that $c_1=c_2=c_3=0, c_4=b_4, c_5=0, c_6=b_6$. Thus by Equation \eqref{eqn09} we have $b_4
+b_6\leq 4$ in particular, which reduces the list to the triples
\begin{align*}
(b_4,b_6,b_8)\in \{(1,0,1),
(2,1,2), (3,0,4)\}
\end{align*}
The first triple yields that $M$ is homothetic to $\hh\pp^4$ due to \cite{GS96}.

For the case $b_4=2$ we may cite theorem \cite[Theorem 1.1, p.~2]{Her06},
which excludes exactly a possible occurrence of this Betti number configuration.

If $b_4=3$, as above we derive that
\begin{align*}
\dim \pi_\textrm{odd}(M)\otimes \qq=\dim
\pi_\textrm{even}(M)\otimes \qq=3
\end{align*}
We compute
\begin{align*}
\dim \big(\Lambda V\big)^8
=&\dim \Sym_2(V^4) +
c_8\\=&\frac{b_4(b_4+1)}{2}+c_8\\=&6+c_8
\end{align*}
As $\dif|_{V^7}$ is injective, we obtain that
\begin{align*}
4=\dim H^8(M)= \dim \big(\Lambda V\big)^8-\dim V^7=6+c_8-c_7
\end{align*}
As $\dim\pi_\textrm{even}(M)\otimes \qq=3$, we see that $c_8=0$ and
that $c_7=2$.  By \eqref{eqn04} we compute that $c_{11}=1$ and
the remaining $c_i$ vanish. This configuration corresponds to $\widetilde \Gr_4(\rr^8)$ as a simple calculation of the minimal model of this space shows.
\end{prf}
A similar investigation for $12$-dimensional positive quaternion K\"ahler manifolds $M$ shows that $M$ has either the Betti numbers and the homotopy Betti numbers of $\hh\pp^3$, $\Gr_2(\cc^5)$, $\widetilde \Gr_4(\rr^7)$ or $b_4=3$, $c_4=3$, $c_7=3$ (with the remaining ones zero).

In fact, from \eqref{eqn02} we derive that $c_4\leq 3$. If $c_4=1$, we know that the cohomology module is concentrated in degrees $0,4,8,12$---we may assume rational $3$-connectedness; and we use $b_6=0$, the fact that odd Betti numbers vanish and Poincar\'e duality for this. 
Due to Hard-Lefschetz it follows that $M\simeq_\qq \hh\pp^3$ in this case.

Assume now that $c_4=2$. Then $b_4=2$ and $b_8=2$ by Poincar\'e duality. Due to Hard-Lefschetz $H^8(M)$ is generated by multiples of the Kraines form. Thus there is a homotopy group in degree $7$ representing a relation in degree $8$, i.e.~$c_7=1$. From \eqref{eqn04} it follows that $c_{15}=1$ and the remaining $c_i$ vanish.

Suppose now that $c_4=3$. Similarly, we derive that $b_4=b_8=3$ and $c_8=0$, $c_7=({3 \choose 2}+3)-3=3$, due to Hard-Lefschetz. It follows from \eqref{eqn05} that $c_i=0$ otherwise.


\section{Geometrically formal manifolds of special holonomy}\label{sec05}

We shall now deal with formal metrics of special holonomy on compact closed manifolds.
Recall that a Riemannian metric is called \emph{formal}, if the product of harmonic forms is again harmonic. A manifold is called \emph{geometrically formal}, if it admits a formal metric. The Betti numbers of geometrically formal manifolds are restricted by the ones of the torus of the same dimension \cite[Theorem 6, p.~524]{Kot01}. In the following we strengthen this under special holonomy. Recall that manifolds of special holonomy are conjectured to be formal (cf.~\cite{AK12}). Formality is a natural obstruction to geometric formality.

\subsection{Manifolds of holonomy $\G_2$ respectively $\Spin(7)$}

In the following we assume the metric of special holonomy to be formal.
\begin{prop}
\begin{itemize}
\item
A compact closed Riemannian manifold with a formal metric of holonomy $\G_2$ satisfies $b_2\leq 14$, $b_3\leq 28$.
\item
A compact closed Riemannian manifold with a formal metric of holonomy $\Spin(7)$ satisfies $b_3\leq 48$, $b_2\leq 21$, $b_4\leq 63$.
\end{itemize}
\end{prop}
\begin{prf}
This is a direct consequence of Theorems \cite[Theorem 10.2.4, p.~246]{Joy00} and \cite[Theorem 10.1.4, p.~244]{Joy00}. Indeed, these theorems yield an orthogonal decomposition of the harmonic forms into subspaces coming from an orthogonal decomposition of the differential forms into irreducible $\G_2$-representations. It is observed that in degree $2$ the subbundle of $\Lambda T^*M$ yielding harmonic forms is at most $14$-dimensional, in dimension three it is at most $28$-dimensional. From the proof of \cite[Theorem 6, p.~524]{Kot01} we know that, since the metric is also formal, this orthogonal decomposition is preserved by harmonic forms and harmonic forms are determined by their coefficients on a fibre. This implies that the dimension of the space of harmonic forms is at most $14$ in degree $2$ and at most $28$ in degree $3$. The Hodge decomposition yields the result.

In the case of $\Spin(7)$ the same arguments together with \cite[Theorem 10.6.5, p.~260]{Joy00} yield the result.
\end{prf}

\subsection{Positive quaternion K\"ahler manifolds and K\"ahler manifolds with trivial Hodge decomposition}\label{sec08}

We call the Hodge decomposition on a K\"ahler manifold trivial, if only $(p,p)$-cohomology exists.
We shall call a Riemannian metric $g$ on a Riemannian manifold $(M^n,g)$ \emph{weakly $p$-formal} for $1\leq p\leq n-1$, if the product of any two harmonic forms in $\ADR(M)^{\leq p}$ is a harmonic form.
We shall call a Riemannian metric $g$ on a Riemannian manifold $(M^n,g)$ \emph{$p$-formal} for $1\leq p\leq n-1$ if the product of a harmonic $i$-form of degree $i\leq p$ with any other harmonic form is harmonic again.

Note that our definition of $p$-formality differs from the one used in \cite{Nag06}.

However, there is a simple and nice parallel to the algebraic concept of classical $s$-formality on compact manifolds---see \cite{FM05} and \cite[Theorem 3.1]{FM05} in particular.
\begin{prop}\label{xprop01}
A Riemannian metric on an ($r$-connected) manifold $M^{n}$ is formal if and only if it is $\lfloor n/2\rfloor$-formal (if and only if it is weakly $(n-r-1)$-formal).
\end{prop}
\begin{prf}
The product of any two harmonic forms is harmonic, since the product of any two forms of degree larger than $\lfloor n/2\rfloor$ vanishes.
\end{prf}

We observe that the proof of the fact that the Betti numbers of a geometrically formal manifold are bounded by the Betti numbers of the torus of the same dimension directly applies to yield a graded version
\begin{prop}\label{xprop02}
If a Riemannian metric on a manifold $M^{n}$ is $p$-formal, then
\begin{align*}
b_k(M)\leq b_k(T^n)
\end{align*}
for $k\leq p$.
\end{prop}
\begin{prf}
Given an $i$-form $x$ with $i\leq p$ we obtain that $x\wedge \ast x=|x|^2 \operatorname{dvol}_g$ is harmonic again and thus $x$ has constant length. Literally reproducing the line of arguments in the proofs of \cite[Lemmas 4 and 5, p.~523]{Kot01} for $i$-forms proves the result in our case.
\end{prf}

 \vspace{5mm}

We generalise the observation made in \cite[Corollary 3.1, p.~208]{Nag06} saying that primitive elements in $H^{1,1}(M)$ vanish if and only if their square vanishes using different techniques. This will improve an estimate on Betti numbers (cf.~\cite[Theorem 6, p.~524]{Kot01} for the general case and \cite[Corollary 4.1, p.~211]{Nag06} for the K\"ahler case).
\begin{theo}\label{xtheo01}
Let $(M^{2n},g)$ be a K\"ahler manifold with Hodge decomposition satisfying $H^*(M;\cc)=\oplus_{i=0}^nH^{i,i}(M)$, i.e.~mixed terms of the Hodge decomposition vanish. Suppose further that the metric $g$ is $p$-formal.

Then the Betti numbers of $M$ satisfy $b_0=b_{2n}=1$, $b_{\odd}=0$ and
\begin{align*}
&b_{2k}(M)=b_{2(n-k)}\leq {n\choose k}\cdot \left({n \choose k} - \bigg(\sum_{i=2}^{\lfloor n/k\rfloor-(\lfloor n/k\rfloor \mod 2)} {n-i+1\choose k-1}\bigg)\right)
\intertext{respectively}
&b_{2k}(M)=b_{2(n-k)}\\&\leq {n\choose k}^2-\sum_{i=2}^{\lfloor n/k\rfloor-(\lfloor n/k\rfloor \mod 2)} {n-\lfloor i/2\rfloor ) \choose k- (i\mod 2)}\cdot{n-\lfloor (i-1)/2\rfloor  \choose k-((i-1) \mod 2))}
\end{align*}
for $0<2k\leq \min\{p, n\}$. In particular, if $(M,g)$ is geometrically formal, this holds for $0<2k\leq n$.

If the metric is $2$-formal, we can improve the estimate on the $2$-forms as
\begin{align*}
b_2(M)\leq n-1
\end{align*}
More generally, if $n$ is divisible by $p$ and if $(M,g)$ is $p$-formal ($p$ even), we obtain
\begin{align*}
b_p(M)\leq {n-1 \choose p/2-1}\cdot {n \choose p/2}
\end{align*}
\end{theo}
\begin{prf}
We shall apply the Hodge--Riemann bilinear relations on a K\"ahler manifold:
From \cite[Theorem V.6.1, p.~203]{Wel08} we cite that---on a compact
K\"ahler manifold $M^{2n}$ (of real dimension $2k$)---the form
\begin{align*}
\tilde Q(\eta,\mu)=\sum_{\max\{s\geq (r-n),0\}} (-1)^{[r(r+1)/2]+s}
\int_M L^{n-r+2s}(\eta_s \wedge \mu_s)
\end{align*}
with Lefschetz decompositions  $\eta=\sum L^s \eta_s\in H^r(M;\cc)$ and
$\mu=\sum L^s\mu_s\in H^r(M;\cc)$, i.e.~with primitive $\eta_s$,
$\mu_s$, satisfies
\begin{align*}
\tilde Q(\eta,J\bar\eta)>0
\end{align*}
for $\eta\neq 0$ and $J=\sum_{a,b} i^{a-b}\Pi_{a,b}$ with canonical
projections \linebreak[4]$\Pi_{a,b}\co H^{a+b}(M;\cc)\to H^{a,b}(M)$ and Lefschetz operator $L$.

This, together with the Lefschetz decomposition, implies that for a K\"ahler manifold with trivial Hodge decomposition---in which case $\tilde Q(\eta, J\bar \eta)=\tilde Q(\eta,\eta)$---the \emph{generalised intersection form} $\tilde Q$ is positive definite.

We shall use this in the following way: The exterior square of a non-trivial $i$-form with $i\leq n$ cannot vanish.

\vspace{5mm}

From \cite[Lemma 5, p.~524]{Kot01} we recall that a form on $M$ is harmonic if and only if it is the linear combination over $C^\infty(M,\rr)$ of an orthonormal system of constant length harmonic forms with \emph{constant} coefficient functions only. We adapt the arguments presented there (requiring the formality of the metric) to our case as in Proposition \ref{xprop02}.

Thus, if $(M,g)$ is $p$-harmonic, the vector space of harmonic $i$-forms with $i\leq p$ is restricted by the number $2n \choose i$ of $i$-forms on $\rr^{2n}$.

The $p$-formality of the metric tells us that arbitrary powers of a harmonic form $x\in H^{2k}(M;\cc)$ are harmonic provided $2k \leq  p$.

Using the fact that $H^{2k}(M;\cc)=H^{k,k}(M)$ we decompose
\begin{align*}
x=\sum_{j} a_j \dif z_{j_1}\wedge \dots \wedge \dif z_{j_k} \wedge \dif \bar z_{j'_1}\wedge \dots \wedge \dif \bar z_{j'_k}
\end{align*}
with $a_j\in \cc$. Assume $x$ to be non-trivial. We derive that $x^{s}$ with
\begin{align}\label{xeqn01}
s&\leq \lfloor n/k\rfloor\\
\nonumber s &\equiv 0 \mod 2
\end{align}
is a non-vanishing harmonic form---the fact that $s$ is even lets us write $x^s$ as a square.

\vspace{5mm}

Let us now investigate what it means for $x$ that a certain power of it does not vanish. Since $(\dif z_i)^2=(\dif \bar z_i)^2=0$, we derive that $x^s\neq 0$ implies that $x$ has at least $s$ non-zero summands (i.e.~in particular $l\geq s$) with the property that they pairwise do not share any coordinate forms $\dif z_i, \dif \bar z_i$ as factors in their respective wedge products. Thus we are looking for a linear subspace $V^{(s)}$ of the $(k,k)$-forms $\ADR(\cc^{k,k})$ with the property that every non-trivial element in there has at least $s$ summands and its $s$-fold power is non-vanishing.

We shall provide an estimate for the dimension of this space. For this we determine a direct linear complement to the described space. We observe that every element in the space
\begin{align*}
C^{(2)}&:=\langle \dif z_1\wedge \dif z_2 \wedge \ldots \wedge \dif \bar z_1 \wedge \dots ,\dif z_1\wedge \dif z_3 \wedge \dots ,  \dots\rangle
\intertext{generated by $(k,k)$-forms which all contain $\dif z_1$ as a factor has a vanishing square.
We set}
C^{(3)}&:=C^{(2)}\oplus \langle \dif z_2\wedge \dif z_3 \wedge \ldots \wedge \dif \bar z_1 \wedge \dots ,\dif z_2\wedge \dif z_4 \wedge \dots ,  \dots\rangle\\
&\vdots
\\
C^{(s)}&:=C^{(s-1)}\\&\oplus\langle \dif z_{s}\wedge \dif z_{s+1} \wedge \ldots \wedge \dif \bar z_1 \wedge \dots ,\dif z_{s}\wedge \dif z_{s+2} \wedge \dots,  \dots\rangle
\end{align*}
That is, the spaces with index $l$ consist of all $(k,k)$-forms containing the factor $\dif z_{l}$ and not containing the ones fixed in the previous factors.
We observe that every element in $C^{(l)}$ has vanishing $l$-fold power. Indeed, let $x_i \in C^{(i)}$. Then $(\sum_i x_{i})^l$ is a linear combination of elements of word-length $kl$ in the $\dif z_i$. Due to the construction of the $C^{(i)}$, for every summand there is an $i$ such that $\dif z_i$ or $\dif \bar z_i$ appears twice as a factor; i.e.~every summand vanishes. In particular, the intersection of $C^{(s)}$ with $V^{(s)}$ must be trivial.

We compute the dimensions of the $C^{(i)}$ as
\begin{align*}
\dim C^{(2)}&={n-1\choose k-1}\cdot{n \choose k}\\
\dim C^{(3)}&=\dim C^{(2)}+{n-2\choose k-1}\cdot{n \choose k}\\
&\vdots\\
\dim C^{(s)}&=\sum_{i=2}^{s-1} \dim C^{(i)}+{n-s+1\choose k-1}\cdot{n \choose k}=\bigg(\sum_{i=2}^{s} {n-i+1\choose k-1}\bigg)\cdot {n\choose k}
\end{align*}
For this, due to Equation \eqref{eqn01}, we may assume that $s+k\leq n$ if $1<k<n$. Let us discuss the cases $k=1$ and $k=n$ separately. In the latter case we only obtain $s=0$ (under the restrictions of the relations \eqref{eqn01})  and $\dim H^{0,0}(M)=1$. If $k=1$, the estimate holds as well.

Thus the dimension of $V^{(s)}$ satisfies
\begin{align*}
\dim V^{(s)}&\leq \dim \ADR(\cc^{k,k}) - \dim C^{(s)}\\&={n\choose k}\cdot \bigg({n \choose k} - \bigg(\sum_{i=2}^{s} {n-i+1\choose k-1}\bigg)\bigg)
\end{align*}
if $1<k<n$.

We assume $k\leq \lfloor n/2\rfloor$ and we let $0<\tilde s$ be the largest $s$ satisfying the relations \eqref{eqn01}; obviously, we have
\begin{align*}
\tilde s=\lfloor n/k\rfloor-(\lfloor n/k\rfloor \mod 2)
\end{align*}
With the vector space $\mathcal{H}^{k,k}(M)$ of harmonic ${k,k}$-forms we derive
\begin{align*}
b_{2k}(M)&\leq \dim \mathcal{H}^{k,k}(M)\\& \leq \dim V^{\tilde s}\\&\leq {n\choose k}\cdot \left({n \choose k} - \bigg(\sum_{i=2}^{\lfloor n/k\rfloor-(\lfloor n/k\rfloor \mod 2)} {n-i+1\choose k-1}\bigg)\right)
\end{align*}
This estimate together with Poincar\'e duality yields the first assertion.

\vspace{5mm}

It is obvious how to produce variations of the spaces $C^{(l)}$. We shall only produce one of these, which will yield the second estimate.
\begin{align*}
C^{(2)}&:=\langle \dif z_1\wedge \dif z_2 \wedge \dots \wedge \dif \bar z_1 \wedge \dots ,\dif z_1\wedge \dif z_3 \wedge \dots ,  \dots\rangle
\intertext{
generated by $(k,k)$-forms which all contain $\dif z_1$ as a factor has a vanishing square.
We set}
C^{(3)}&:=C^{(2)}\oplus \langle \dif \bar z_1\wedge \dif z_2 \wedge \dots ,\dif \bar z_1\wedge \dif z_3 \wedge \dots ,  \dots\rangle\\
C^{(4)}&:=C^{(3)}\oplus \langle \dif z_2\wedge \dif z_3 \wedge \dots ,\dif  z_2\wedge \dif z_4 \wedge \dots ,  \dots\rangle\\
C^{(5)}&:=C^{(4)}\oplus \langle \dif \bar z_2\wedge \dif z_3 \wedge \dots ,\dif \bar z_2\wedge \dif z_4 \wedge \dots ,  \dots\rangle\\
&\vdots
\\
C^{(s)}&:=C^{(s-1)}\\&\oplus\langle \dif z_{s/2}\wedge \dif z_{s/2+1} \wedge \dots \wedge \dif \bar z_1 \wedge \dots ,\dif z_{s/2}\wedge \dif z_{s/2+2} \wedge \dots,  \dots\rangle
\end{align*}
(if $s$ is even.) That is, the spaces with even index $l$ consist of all $(k,k)$-forms containing the factor $\dif z_{l/2}$ and not containing the ones fixed in the previous factors. The spaces with odd index consist of all forms fixing the element $\dif \bar z_{(l-1)/2}$ and not containing any of the previous forms.
Again we observe that every element in $C^{(l)}$ has vanishing $l$-fold power.

We compute the dimensions of the $C^{(i)}$ as
\begin{align*}
\dim C^{(2)}&={n-1\choose k-1}\cdot{n \choose k}\\
\dim C^{(3)}&=\dim C^{(2)}+{n-1\choose k-1}\cdot{n-1 \choose k}\\
\dim C^{(4)}&=\dim C^{(3)}+{n-2\choose k}\cdot{n-1 \choose k-1}\\
\dim C^{(5)}&=\dim C^{(4)}+{n-2\choose k-1}\cdot{n-2 \choose k}\\
&\vdots\\
\dim C^{(s)}&=\sum_{i=2}^{s-1} \dim C^{(i)}+{n-\lfloor s/2\rfloor ) \choose k- (s\mod 2)}\cdot{n-\lfloor (s-1)/2\rfloor  \choose k-((s-1) \mod 2))}
\\&=\sum_{i=2}^{s} {n-\lfloor i/2\rfloor ) \choose k- (i\mod 2)}\cdot{n-\lfloor (i-1)/2\rfloor  \choose k-((i-1) \mod 2))}
\end{align*}

Thus the dimension of $V^{(s)}$ satisfies
\begin{align*}
\dim V^{(s)}&\leq \dim \ADR(\cc^{k,k}) - \dim C^{(s)}\\&={n\choose k}^2-\sum_{i=2}^{s} {n-\lfloor i/2\rfloor ) \choose k- (i\mod 2)}\cdot{n-\lfloor (i-1)/2\rfloor  \choose k-((i-1) \mod 2))}
\end{align*}
if $1<k<n$ and the asserted second estimate follows.

\vspace{5mm}

Let us improve this estimate for $b_2(M)$ under the assumption of $2$-formality and $n$ being even. Let $x=\sum_{j}^k a_j \dif z_{j_1}\wedge \dif \bar z_{j_2}\in \mathcal{H}^{1,1}(M)$ be a harmonic $(1,1)$-form with (constant) complex coefficients $a_j\in \cc$. Then $x^n$ is a non-vanishing harmonic $2n$-form. Thus $x$ has at least $n$ non-vanishing summands with the (in total $2n$) pairwise disjoint factors $\dif z_i, \dif\bar z_i$.
That is, every such form $x$, i.e.~every sum of monomials, contains a subsum uniquely determined by two permutations $\pi, \tau$ of $\{1,2,\dots,n\}$. In other words, $x=a_1\cdot \dif z_{\pi(1)}\wedge \dif \bar z_{\tau(1)} + a_2\cdot \dif z_{\pi(2)}\wedge \dif \bar z_{\tau(2)} + \dots + a_n\cdot \dif z_{\pi(n)}\wedge \dif \bar z_{\tau(n)}+ \dots$ (and other summands).

We shall now find an upper bound for the dimension of a vector space $V$ containing only elements of the form of $x$. We shall see that $\dim  V\leq n-1$. Indeed, without restriction, we may choose $n$ basis elements $(x_i)_{1\leq i\leq n}$ of the above form. Expressing these elements in the basis $(\dif z_i\wedge \dif \bar z_j)_{1\leq i,j\leq n}$ yields an $n\times {n\choose 2}$ matrix. The first $n$ columns give the coefficients of the elements $(\dif z_1\wedge \dif \bar z_j)_{1\leq j\leq n}$. (No further $\dif z_1$ appear.) We may bring this matrix into ``upper triangular form''---by abuse of notation, as it not quadratic. This implies that we may find a linear combination $x$ of the $x_i$ with the property that $x$ has no non-trivial summand containing $\dif z_1$ as a factor. Consequently, $x^n=0$; a contradiction.

The analogous arguments apply for the respective estimate on $p$-forms. We use that after fixing $\dif z_1$ we may choose $p/2-1$ coordinate funtions out of $\dif z_i$ and $p/2$-many out of the $\dif \bar z_i$.

\vspace{5mm}

Finally, under the assumption of $2$-formality let us give the line of arguments which provides that $b_2(M)\leq n-1$ irrespective of the parity of $n$. For this we consider the Riemannian product $M\times \cc\pp^1$ which is of dimension $2(n+1)$. It is easy to see---cf.~\cite[Lemma A.2, p.~360]{OP11}---that its harmonic $(1,1)$-forms are given by the linear combinations of the (pullbacks of the Euclidean projections of the) harmonic $2$-forms of $M$ and $\cc\pp^1$. Thus $b_i(M\times \cc\pp^1)=b_i(M)$ for $i\neq 2$, $b_2(M\times \cc\pp^1)=b_2(M)+1$ and $M\times \cc\pp^1$ has trivial Hodge decomposition in the sense of this theorem. If $n$ is odd, we have that $b_2(M\times \cc\pp^1)\leq n$ by the arguments above. This implies that $b_2(M)\leq n-1$. This finishes the proof.
\end{prf}

\begin{cor}
The estimates from Theorem \ref{xtheo01} hold for the twistor space $(Z,g)$ of a positive quaternion K\"ahler manifold provided $(Z,g)$ is $p$-formal.
\end{cor}
\begin{prf}
According to the proof of \cite[Theorem 6.6, p.~163]{Sal82} the K\"ahler manifold $(Z,g)$ has a Hodge decomposition concentrated in $(k,k)$-degrees only.
\end{prf}
The simplest example of such a manifold is $\cc\pp^{2n+1}$. It was shown in \cite{KT09} that the twistor space of $\Gr_2(\cc^{n+2})$ is not geometrically formal in general.

\vspace{5mm}

We shall now state the analogous theorem for positive quaternion K\"ahler manifolds, i.e.~connected oriented manifolds with complete metric, positive scalar curvature and holonomy contained in $\Sp(n)\Sp(1)$---see \cite{Sal82}, \cite{Sal99}, \cite{Ama09}.
\begin{theo}\label{xtheo02}
Let $(M^{4n},g)$ be a positive quaternion K\"ahler manifold. Suppose further that the metric $g$ is $p$-formal.

Then the Betti numbers of $M$ satisfy $b_0=b_{2n}=1$, $b_2\leq 2$, $b_{\odd}=0$ and
\begin{align*}
b_{k}(M)=b_{4n-k}\leq {4n\choose k} - \sum_{i=2}^{\lfloor 4n/k\rfloor-(\lfloor 4n/k\rfloor \mod 2)} {4n-i+1\choose k-1}
\end{align*}
for $0<k\leq \min\{p, 2n\}$. In particular, if $(M,g)$ is geometrically formal, this holds for $0<k\leq 2n$.

If $2n$ is divisible by $p$ and if $(M,g)$ is $p$-formal, we obtain
\begin{align*}
b_p(M)\leq {4n-1 \choose p-1}
\end{align*}
\end{theo}
\begin{prf}
Basically due to the fact that the twistor space of a positive quaternion K\"ahler manifold is a K\"ahler manifold with trivial Hodge decomposition in the above sense, the twistor transform yields that the generalised intersection form
\begin{align*}
Q(x,y)=(-1)^{r/2} \int_M x\wedge y \wedge u^{n-r/2}
\end{align*}
of $M$ (with $[x],[y] \in H^r(M^{4n},\rr)$, $r\geq 0$ even and where $u$ denotes the Kraines form in degree $4$)  is positive definite---see \cite{Fuj87}, \cite{NT83}, \cite[Theorem 1.16, p.~13]{Ama09}. Consequently, as in the proof of Theorem \ref{xtheo01} we obtain that $p$-formality implies that the $s$-fold power $x^s$ of a $k$-form $0\neq x$ with $k\leq p$ is a non-vanishing harmonic form provided
\begin{align}\label{xeqn02}
s&\leq \lfloor 4n/k \rfloor
\intertext{and}
s&\equiv 0 \mod 2
\end{align}

We may reproduce the arguments from the proof of Theorem \ref{xtheo01} in order to show that the spaces
\begin{align*}
C^{(2)}&:=\langle {\dif x_1\wedge \dif x_2\dots,\dif x_1\wedge \dif x_3\wedge \dots, \dots} \rangle\\
&\vdots\\
C^{(s)}&:=C^{(s-1)}\oplus \langle {\dif x_s\wedge \dif x_{s+1}\dots,\dif x_s\wedge \dif x_{s+2}\wedge \dots, \dots} \rangle
\end{align*}
(defined as above) with coordinate forms $\dif x_i$ satsify that an element from $C^{(l)}$ has vanishing $l$-fold power. Again we compute the dimensions and obtain
\begin{align*}
\dim C^{(2)}&={4n-1\choose k-1}\\
\dim C^{(3)}&=\dim C^{(2)}+{4n-2\choose k-1}\\
&\vdots\\
\dim C^{(s)}&=\sum_{i=2}^{s-1} \dim C^{(i)}+{4n-s+1\choose k-1}=\sum_{i=2}^{s} {4n-i+1\choose k-1}
\end{align*}

Define $V^{(s)}$ as a linear subspace of the $k$-forms $\ADR(\cc^{k,k})$ with the property that every non-trivial element in there has at least $s$ summands and its $s$-fold power is non-vanishing. Thus the dimension of $V^{(s)}$ satisfies
\begin{align*}
\dim V^{(s)}&\leq \dim \ADR(\cc^{k,k}) - \dim C^{(s)}\\&={4n\choose k} - \sum_{i=2}^{s} {4n-i+1\choose k-1}
\end{align*}

We assume $k\leq 2n$ and we let $0<\tilde s$ be the largest $s$ satisfying the relations \eqref{eqn02}, i.e.
\begin{align*}
\tilde s=\lfloor 4n/k\rfloor-(\lfloor 4n/k\rfloor \mod 2)
\end{align*}
With the vector space $\mathcal{H}^k(M)$ of harmonic $k$-forms we derive
\begin{align*}
b_{k}(M)&\leq \dim \mathcal{H}^k(M)\\& \leq \dim V^{\tilde s}\\&={4n\choose k} - \sum_{i=2}^{\lfloor 4n/k\rfloor-(\lfloor 4n/k\rfloor \mod 2)} {4n-i+1\choose k-1}
\end{align*}
This estimate together with Poincar\'e duality, the fact that $b_2\leq 1$ and $b_\odd=0$ on a positive quaternion K\"ahler manifold (see \cite{Sal82}, \cite{Sal99}, \cite[Theorem 1.13, p.~11]{Ama09}) yields the first assertion.

As for the remaining assertion, once again, we argue in the same way as in the proof of Theorem \ref{xtheo01}.
\end{prf}
Examples of geometrically formal positive quaternion K\"ahler manifolds are the symmetric ones, i.e.~the Wolf spaces like $\hh\pp^n$, $\Gr_2(\cc^{n+2})$, $\widetilde \Gr_4(\rr^{n+4})$, etc.

\begin{rem}\label{xrem01}
First of all one needs to observe that building different spaces in the proof of Theorem \ref{xtheo01} will yield several more estimates out of which we only selected two extremal cases. Thus there is still potential for improving the bounds in several specific cases. It is, however, not our goal to completely finetune the results, but to present the ideas behind such a process which is left to the interested reader.

Note that the first estimate in Theorem \ref{xtheo01} yields
\begin{align*}
b_2(M^{2n})&\leq \begin{cases} n & \textrm{if $n$ is even}\\
2n & \textrm{if $n$ is odd}
\end{cases}
\intertext{whereas the second one produces}
b_2(M^{2n})&\leq n^2- \sum_{i=2}^{n-(n \mod 2)} n- \lfloor i/2\rfloor
\intertext{The one in Theorem \ref{xtheo02} yields}
b_2(M^{4n})&\leq 2n^2+n
\end{align*}

In Table \ref{xtable01} we shall explicitly compute some estimates provided by the Theorems---thus comparing the new estimates to the known one provided by $b_i(M^n)\leq b_i(T^{\dim M})$.
\begin{table}[h]
\centering \caption{Exemplary Betti number bounds} \label{xtable01}
\begin{tabular}{@{\hspace{3mm}}c@{\hspace{3mm}}| @{\hspace{3mm}}l@{\hspace{3mm}}|@{\hspace{3mm}}l@{\hspace{3mm}}|@{\hspace{3mm}}l}
$\dim M$= &  K\"ahler trivial Hodge& PQK & $T^{\dim M}$\\
\hline
$4$ & $b_2\leq 3 $ & & $b_2=6 $\\
$6$ & $b_2\leq 5 $ & & $b_2=15 $\\
$8$ & $b_2\leq 7 $ & & $b_2=28 $\\
& $b_4\leq 18 $ &  & $b_4= 70$\\
$10$ & $b_2\leq 9 $ & & $b_2= 45$\\
& $b_4\leq 60 $ &  & $b_4= 210$\\
$12$ & $b_2\leq 11 $ & $b_2 \leq 2$& $b_2= 66$\\
& $b_4\leq 150$ &  $b_4\leq 165$& $b_4=495 $ \\
& $b_6\leq 200$ &  $b_6\leq 462$ & $b_6=924 $\\
$14$ & $b_2\leq 13$ & & $b_2=91 $\\
& $b_4\leq 315$ &  & $b_4= 1001$\\
& $b_6\leq 700$ &   & $b_6= 3003$\\
$16$ & $b_2\leq 15 $ & $b_2\leq 2$& $b_2= 120$\\
& $b_4\leq 196$ &  $b_4\leq 455$& $b_4=1820 $ \\
& $b_6\leq 1960 $ &  $b_6\leq 5005 $ & $b_6= 8008$\\
& $b_8\leq 2450 $ &  $b_8\leq 6435$ & $b_8= 12870$\\
\end{tabular}
\end{table}
Note that due to the classification of positive quaternion K\"ahler manifolds in dimensions $4$ and $8$, we start computing bounds for them in dimension $12$.
Moreover, one might now apply a known relation on the Betti numbers of a positive quaternion K\"ahler manifolds---\cite{Sal82}, \cite{Sal99}, \cite[Theorem 1.13, p.~11]{Ama09}---in order to improve the situation slightly further.
\end{rem}

\vspace{3mm}

Let us now motivate the conjecture we made in the introduction and which claims that positive quaternion K\"ahler manifolds are both (rationally) elliptic and geometrically formal.
Indeed, the LeBrun--Salamon conjecture implies that a positive quaternion K\"ahler manifold should be a symmetric space. Consequently, it is homogeneous and elliptic and geometrically formal, since the invariant forms are the harmonic forms due to Cartan.

There is also direct instrinsic motivation for this conjecture. Note that the curvature condition is actually stronger than just positive scalar curvature: The sum of the curvatures in $I$, $J$, $K$ directions are positive (\cite[Formula 14.42b, p.~406]{Bes08}). Moreover, important results from the theory of positive sectional curvature like the \emph{connectedness theorem} by Wilking (\cite{Wil03}) also hold for Positive Quaternion K\"ahler manifolds (\cite[Theorem A, p.~150]{Fan04}). This might suggest that in the quaternionic setting positive
scalar curvature might be regarded as a substitute for positive sectional curvature to a certain extent.

The \emph{Bott conjecture} claims that a manifold of non-negative sectional curvature  should be (rationally) elliptic (\cite{Gro93}). The equally famous \emph{Hopf conjecture} in positive curvature states that the Euler characteristic should be positive. A combination of both says that the manifold should be positively elliptic. In the quaternionic setting, we already have positive Euler characteristic; the fact that the odd Betti numbers of a positive quaternion K\"ahler manifolds vanish would be just equivalent to positive Euler characteristic once we have rational ellipticity (see \cite[Proposition 32.10, p.~444]{FHT01}). In other words, this obstruction to rational ellipticity vanishes. The first part of the conjecture above hence can be considered a \emph{quaternionic Bott conjecture}. (Of course, note that if a quaternion K\"ahler manifold has positive sectional curvature, it is homothetic to $\hh\pp^n$ already.)

It seems that until now not many more examples of geometrically formal manifolds are known than symmetric spaces, some homogeneous ones, etc (cf.~\cite[Theorem 13, p.~503]{KT03}). A confirmation of the conjecture above then might be at least morally a good progress towards the LeBrun--Salamon conjecture. (Note that the homogeneous positive quaternion K\"ahler manifolds were classified to be Wolf spaces by Alekseevski.)

We observe that both implications of the conjecture separately imply upper bounds on the Betti numbers, as we have seen above. Till today no general upper bound even just on the Euler characteristic seems to be known.

However, not only the implications are similar, both claims also share a common obstruction, which is formality. Indeed, it is easy to see that formality is an obstruction to geometric formality (using the Hodge decomposition). Due to Halperin it is know that rationally elliptic spaces of positive Euler characteristic are formal. This obstruction of formality vanishes on arbitrary positive quaternion K\"ahler manifolds due to \cite{AK12}.

\section{Changing coefficients}\label{sec09}

Let $\Bbbk$ and $\kk\Ni \Bbbk$ be fields of characteristic $0$. Denote by
\begin{align*}
\mathfrak{m}_{\Bbbk\In \kk}(X):=\{\textrm{$\Bbbk$-homotopy types which have the $\kk$-homotopy type $X_{(\kk)}$} \}
\end{align*}
for a simply-connected CW-complex $X$. (By $X_{(\kk)}$ we denote localisation, as usual---see \cite[Chapter 9, p.~102]{FHT01}.)

\begin{theo}\label{theo05}
Suppose that $\kk\Ni \Bbbk$ is a finite field extension. Assume further that $\dim H^*(X;\qq)$ is finite dimensional. 

Then $|\mathfrak{m}_{\Bbbk\In \kk}(X)|$ is finite. In particular, there are finitely many $n$-dimensional $\Bbbk$-types if and only if there are finitely many such $\kk$-types.
\end{theo}
\begin{prf}
Let $(\Lambda V,\dif)$ be the minimal Sullivan model of $X$ over $\Bbbk$. Then $(\Lambda V\otimes \kk, \dif\otimes \kk)$ is the minimal Sullivan model of $X$ over $\kk$.

Let us show that the number of $\Bbbk$-automorphisms of $(\Lambda V,\dif)$ tensoring to the same $\kk$-automorphism is finite. This will imply the result.

For this write $\kk=\Bbbk[r_1,r_2,\dots, r_l]$---$\kk$ is a finite field extension of $\Bbbk$.
Denote by $f_{x_i,j}$ the morphism $\Lambda V\otimes \kk\to \Lambda V\otimes \kk$, $x_i\mapsto r_j x_i$ for $x_i$, $i\in I$, a homogeneous basis element of $V$ and $1\leq j\leq l$.
Consider the monoid
\begin{align*}
A:=(f_{x_i,j} \mid i\in I, 1\leq j\leq l )\cup \{\id\}
\end{align*}
generated via composition by the $f_{x_i,j}$. Then every automorphism of $\Lambda V\otimes \kk$ is of the form $g\circ (f\otimes \kk)$ with $f\in \Aut(\Lambda V)$ and $g\in A$.

Let us now understand the automorphisms compatible with differentials. Suppose that $f\in \Aut(\Lambda V,\dif)$ is an automorphism compatible with the differential. Since $H(\Lambda V,\dif)=H^*(X;\qq)\otimes \Bbbk$ is finite-dimensional, we derive that $f(v)$ for $v\in V^{> n}$ is uniquely determined by $f|_{(\Lambda V,\dif)^{\leq n}}$ for $n=\dim X$ the formal dimension of $X$. This implies that an automorphism of $(\Lambda V,\dif)$ is determined by a suitable choice of the image of a finite homogeneous basis of $V^{\leq n}$. From the fact that the automorphisms of $\Lambda V\otimes \kk$ decompose as $g\circ (f\otimes \kk)$ with $g\in A$, $f\in \Aut(\Lambda V)$, we deduce that there are only finitely many automorphisms of $(\Lambda V\otimes \kk,\dif\otimes \kk)$ which are not of the form $h\otimes \kk$ with $h\in \Aut(\Lambda V,\dif)$.

Every rational type $(\Lambda \tilde V,\tilde\dif)\in \mathfrak{m}_{\Bbbk\In \kk}(X)$ induces an isomorphism $(\Lambda \tilde V\otimes \kk,\tilde\dif\otimes \kk)\cong (\Lambda V\otimes \kk,\dif\otimes \kk)$ of $\kk$-types, respectively an automorphism of $(\Lambda V\otimes \kk,\dif\otimes \kk)$. If this automorphism is induced by an automorphism over $\Bbbk$, then the two types are actually one type over $\Bbbk$. Since there are only finitely many such automorphisms, there are only finitely many $\Bbbk$-types, i.e.~$\mathfrak{m}_{\Bbbk\In \kk}(X)$ is finite.

The assertion is trivial, if there are only finitely many $\kk$-types.
\end{prf}



\begin{thebibliography}{10}

\bibitem{Ama09}
M.~Amann.
\newblock {\em Positive Quaternion K\"ahler Manifolds}.
\newblock PhD thesis, WWU M\"unster, 2009.
\newblock http://miami.uni-muenster.de/servlets/DocumentServlet?id=4869.

\bibitem{Ama11a}
M.~Amann.
\newblock Positive quaternion {K}\"ahler manifolds with fourth {B}etti number
  equal to one.
\newblock {\em Topology Appl.}, 158(2):183--189, 2011.

\bibitem{Ama12}
M.~Amann.
\newblock Partial classification results for positive quaternion {K}\"ahler
  manifolds.
\newblock {\em International Journal of Mathematics}, 23(2):1--39, 2012.

\bibitem{AD10}
M.~Amann and A.~Dessai.
\newblock The {$\hat A$}-genus of {$S^1$}-manifolds with finite second homotopy
  group.
\newblock {\em C. R. Math. Acad. Sci. Paris}, 348(5-6):283--285, 2010.

\bibitem{AK12}
M.~Amann and V.~Kapovitch.
\newblock On fibrations with formal elliptic fibers.
\newblock {\em Adv. Math.}, 231(3-4):2048--2068, 2012.

\bibitem{AH70}
M.~Atiyah and F.~Hirzebruch.
\newblock Spin-manifolds and group actions.
\newblock In {\em Essays on {T}opology and {R}elated {T}opics ({M}\'emoires
  d\'edi\'es \`a {G}eorges de {R}ham)}, pages 18--28. Springer, New York, 1970.

\bibitem{Bes08}
A.~L. Besse.
\newblock {\em Einstein manifolds}.
\newblock Classics in Mathematics. Springer-Verlag, Berlin, 2008.
\newblock Reprint of the 1987 edition.

\bibitem{BD49}
A.~Borel and J.~De~Siebenthal.
\newblock Les sous-groupes ferm\'es de rang maximum des groupes de {L}ie clos.
\newblock {\em Comment. Math. Helv.}, 23:200--221, 1949.

\bibitem{Cav05}
G.~R. Cavalcanti.
\newblock {\em New aspects of the $\dif \dif_c$-lemma}.
\newblock PhD thesis, Oxford University, 2005.

\bibitem{DGMS75}
P.~Deligne, P.~Griffiths, J.~Morgan, and D.~Sullivan.
\newblock Real homotopy theory of {K}\"ahler manifolds.
\newblock {\em Invent. Math.}, 29(3):245--274, 1975.

\bibitem{Dev11}
J.~DeVito.
\newblock {\em The classification of simply connected biquotients of dimension
  at most 7 and 3 new examples of almost positively curved manifolds}.
\newblock PhD thesis, University of Pennsylvania, 2011.
\newblock http://repository.upenn.edu/edissertations/311/.

\bibitem{DK90}
S.~K. Donaldson and P.~B. Kronheimer.
\newblock {\em The geometry of four-manifolds}.
\newblock Oxford Mathematical Monographs. The Clarendon Press Oxford University
  Press, New York, 1990.
\newblock Oxford Science Publications.

\bibitem{Fan04}
F.~Fang.
\newblock Positive quaternionic {K}\"ahler manifolds and symmetry rank.
\newblock {\em J. Reine Angew. Math.}, 576:149--165, 2004.

\bibitem{Fan08}
F.~Fang.
\newblock Positive quaternionic {K}\"ahler manifolds and symmetry rank. {II}.
\newblock {\em Math. Res. Lett.}, 15(4):641--651, 2008.

\bibitem{FHT91}
Y.~F{\'e}lix, S.~Halperin, and J.-C. Thomas.
\newblock Elliptic spaces.
\newblock {\em Bull. Amer. Math. Soc. (N.S.)}, 25(1):69--73, 1991.

\bibitem{FHT01}
Y.~F{\'e}lix, S.~Halperin, and J.-C. Thomas.
\newblock {\em Rational homotopy theory}, volume 205 of {\em Graduate Texts in
  Mathematics}.
\newblock Springer-Verlag, New York, 2001.

\bibitem{FOT08}
Y.~F{\'e}lix, J.~Oprea, and D.~Tanr{\'e}.
\newblock {\em Algebraic models in geometry}, volume~17 of {\em Oxford Graduate
  Texts in Mathematics}.
\newblock Oxford University Press, Oxford, 2008.

\bibitem{FM05}
M.~Fern{\'a}ndez and V.~Mu{\~n}oz.
\newblock Formality of {D}onaldson submanifolds.
\newblock {\em Math. Z.}, 250(1):149--175, 2005.

\bibitem{Fuj87}
A.~Fujiki.
\newblock On the de {R}ham cohomology group of a compact {K}\"ahler symplectic
  manifold.
\newblock In {\em Algebraic geometry, {S}endai, 1985}, volume~10 of {\em Adv.
  Stud. Pure Math.}, pages 105--165. North-Holland, Amsterdam, 1987.

\bibitem{GS96}
K.~Galicki and S.~Salamon.
\newblock Betti numbers of {$3$}-{S}asakian manifolds.
\newblock {\em Geom. Dedicata}, 63(1):45--68, 1996.

\bibitem{Gro93}
K.~Grove.
\newblock Critical point theory for distance functions.
\newblock In {\em Differential geometry: {R}iemannian geometry ({L}os
  {A}ngeles, {CA}, 1990)}, volume~54 of {\em Proc. Sympos. Pure Math.}, pages
  357--385. Amer. Math. Soc., Providence, RI, 1993.

\bibitem{Her06}
R.~Herrera.
\newblock Positive quaternion-{K}\"ahler 16-manifolds with {$b\sb 2=0$}.
\newblock {\em Q. J. Math.}, 57(2):203--214, 2006.

\bibitem{Joy00}
D.~D. Joyce.
\newblock {\em Compact manifolds with special holonomy}.
\newblock Oxford Mathematical Monographs. Oxford University Press, Oxford,
  2000.

\bibitem{KZ04}
V.~Kapovitch and W.~Ziller.
\newblock Biquotients with singly generated rational cohomology.
\newblock {\em Geom. Dedicata}, 104:149--160, 2004.

\bibitem{Kla88}
S.~Klaus.
\newblock Einfach--zusammenh\"angende kompakte homogene r\"aume bis zur
  dimension neun.
\newblock Master's thesis, Johannes Gutenberg Universit\"at Mainz, 1988.

\bibitem{Kot01}
D.~Kotschick.
\newblock On products of harmonic forms.
\newblock {\em Duke Math. J.}, 107(3):521--531, 2001.

\bibitem{KT03}
D.~Kotschick and S.~Terzi{\'c}.
\newblock On formality of generalized symmetric spaces.
\newblock {\em Math. Proc. Cambridge Philos. Soc.}, 134(3):491--505, 2003.

\bibitem{KT09}
D.~Kotschick and S.~Terzi{\'c}.
\newblock Chern numbers and the geometry of partial flag manifolds.
\newblock {\em Comment. Math. Helv.}, 84(3):587--616, 2009.

\bibitem{NT83}
T.~Nagano and M.~Takeuchi.
\newblock Signature of quaternionic {K}aehler manifolds.
\newblock {\em Proc. Japan Acad. Ser. A Math. Sci.}, 59(8):384--386, 1983.

\bibitem{Nag06}
P.-A. Nagy.
\newblock On length and product of harmonic forms in {K}\"ahler geometry.
\newblock {\em Math. Z.}, 254(1):199--218, 2006.

\bibitem{OP11}
L.~Ornea and M.~Pilca.
\newblock Remarks on the product of harmonic forms.
\newblock {\em Pacific J. Math.}, 250(2):353--363, 2011.

\bibitem{PP03}
G.~P. Paternain and J.~Petean.
\newblock Minimal entropy and collapsing with curvature bounded from below.
\newblock {\em Invent. Math.}, 151(2):415--450, 2003.

\bibitem{Sal82}
S.~Salamon.
\newblock Quaternionic {K}\"ahler manifolds.
\newblock {\em Invent. Math.}, 67(1):143--171, 1982.

\bibitem{Sal89}
S.~Salamon.
\newblock {\em Riemannian geometry and holonomy groups}, volume 201 of {\em
  Pitman Research Notes in Mathematics Series}.
\newblock Longman Scientific \& Technical, Harlow, 1989.

\bibitem{Sal93}
S.~Salamon.
\newblock Index theory and quaternionic {K}\"ahler manifolds.
\newblock In {\em Differential geometry and its applications ({O}pava, 1992)},
  volume~1 of {\em Math. Publ.}, pages 387--404. Silesian Univ. Opava, Opava,
  1993.

\bibitem{Sal99}
S.~Salamon.
\newblock Quaternion-{K}\"ahler geometry.
\newblock In {\em Surveys in differential geometry: essays on {E}instein
  manifolds}, Surv. Differ. Geom., VI, pages 83--121. Int. Press, Boston, MA,
  1999.

\bibitem{Sin93}
W.~Singhof.
\newblock On the topology of double coset manifolds.
\newblock {\em Math. Ann.}, 297(1):133--146, 1993.

\bibitem{Sul77}
D.~Sullivan.
\newblock Infinitesimal computations in topology.
\newblock {\em Inst. Hautes \'Etudes Sci. Publ. Math.}, (47):269--331 (1978),
  1977.

\bibitem{Tot02}
B.~Totaro.
\newblock Cheeger manifolds and the classification of biquotients.
\newblock {\em J. Differential Geom.}, 61(3):397--451, 2002.

\bibitem{Tot03}
B.~Totaro.
\newblock Curvature, diameter, and quotient manifolds.
\newblock {\em Math. Res. Lett.}, 10(2-3):191--203, 2003.

\bibitem{Wel08}
R.~O. Wells, Jr.
\newblock {\em Differential analysis on complex manifolds}, volume~65 of {\em
  Graduate Texts in Mathematics}.
\newblock Springer, New York, third edition, 2008.
\newblock With a new appendix by Oscar Garcia-Prada.

\bibitem{Wil03}
B.~Wilking.
\newblock Torus actions on manifolds of positive sectional curvature.
\newblock {\em Acta Math.}, 191(2):259--297, 2003.

\end{thebibliography}


\vfill

\begin{center}
\noindent
\begin{minipage}{\linewidth}
\small \noindent \textsc
{Manuel Amann} \\
\textsc{Fakult\"at f\"ur Mathematik}\\
\textsc{Institut f\"ur Algebra und Geometrie}\\
\textsc{Karlsruher Institut f\"ur Technologie}\\
\textsc{Kaiserstra\ss e 89--93}\\
\textsc{76133 Karlsruhe}\\
\textsc{Germany}\\
[1ex]
\textsf{manuel.amann@kit.edu}\\
\textsf{http://topology.math.kit.edu/$21\_54$.php}
\end{minipage}
\end{center}

\end{document}